\documentclass{article}
\usepackage[utf8]{inputenc}
\usepackage[left=1in,top=1in,right=1in,bottom=1in]{geometry}
\usepackage[linktocpage=true]{hyperref}
\usepackage{setspace}
\usepackage{amssymb, amsmath, amsthm, graphicx,mathrsfs}
\usepackage{caption,cite}
\usepackage{mathtools,color}

\newtheorem{thm}{Theorem}[section]

\newtheorem{cor}[thm]{Corollary}

\newtheorem{lem}[thm]{Lemma}
\newtheorem{conj}[thm]{Conjecture}
\newtheorem{prop}[thm]{Proposition}
\newtheorem{claim}[thm]{Claim}

\renewcommand{\S}{\mathcal{S}}
\newcommand{\C}{\mathcal{C}}
\newcommand{\F}{\mathcal{F}}
\newcommand{\Fb}{\mathrm{Forb}}

\newcommand{\N}{\mathbb{N}}
\newcommand{\E}{\mathbb{E}}
\newcommand{\B}{\mathcal{B}}

\renewcommand{\P}{\mathrm{Pr}}
\renewcommand{\l}{\left}
\renewcommand{\r}{\right}
\newcommand{\ex}{\mathrm{ex}}

\usepackage[affil-it]{authblk}

\title{Tur\'an theorems for even cycles in random hypergraph}
\author{Jiaxi Nie\thanks{jiaxi\_nie@fudan.edu.cn}}
\affil{Shanghai Center for Mathematical Sciences, Fudan University, Shanghai, China}

\begin{document}

\maketitle

\begin{abstract}
    Let $\F$ be a family of $r$-uniform hypergraphs. The random Tur\'an number $\mathrm{ex}(G^r_{n,p},\F)$ is the maximum number of edges in an $\F$-free subgraph of $G^r_{n,p}$, where $G^r_{n,p}$ is the Erd\H{o}s-R\'enyi random $r$-graph with parameter $p$. Let $C^r_{\ell}$ denote the $r$-uniform linear cycle of length $\ell$. For $p\ge n^{-r+2+o(1)}$, Mubayi and Yepremyan showed that $\mathrm{ex}(G^r_{n,p},C^r_{2\ell})\le\max\{p^{\frac{1}{2\ell-1}}n^{1+\frac{r-1}{2\ell-1}+o(1)},pn^{r-1+o(1)}\}$. This upper bound is not tight when $p\le n^{-r+2+\frac{1}{2\ell-2}+o(1)}$. In this paper, we close the gap for $r\ge 4$. More precisely, we show that $\mathrm{ex}(G^r_{n,p},C^r_{2\ell})=\Theta(pn^{r-1})$ when $p\ge n^{-r+2+\frac{1}{2\ell-1}+o(1)}$. Similar results have recently been obtained independently in a different way by Mubayi and Yepremyan. For $r=3$, we significantly improve Mubayi and Yepremyan's upper bound. Moreover, we give reasonably good upper bounds for the random Tur\'an numbers of Berge even cycles, which improve previous results of Spiro and Verstra\"ete.
\end{abstract}

\section{Introduction}
An {\em $r$-uniform hypergraph}, or {\em $r$-graph} for short, is a hypergraph whose edges all have size $r$. Given a family $\F$ of $r$-graphs, the \emph{Turán number} $\ex(n,{\F})$  is the largest integer $m$ such that there exists an $r$-graph with $n$ vertices and $m$ edges that does not contain any member in $\F$ as a subgraph. The {\em Tur\'an problem} consists of determining or estimating $\ex(n,{\F})$, which is a central problem in extremal combinatorics. Turán in his seminal paper\cite{turan1941} showed that $\ex(n,K_t)=(1-\frac{1}{t-1}+o(1))\binom{n}{2}$, where $K_t$ is the complete graph on $t$ vertices. This result was later generalized by Erd\H{o}s, Stone and Simmonovits\cite{erdos1946structure,Erdos1966ALT} who showed that $\ex(n,G)=(1-\frac{1}{\chi(G)-1}+o(1))\binom{n}{2}$, where $\chi(G)$ is the chromatic number of $G$. This result essentially solved the Tur\'an problems for all graphs with chromatic number at least three. For bipartite graphs, only sporadic results are known. For the graph even cycle, $C_{2\ell}$, Erd\H{o}s~\cite{furedi2013history} and independently Bondy and Simonovits~\cite{bondy1974cycles} showed that $\ex(n, C_{2\ell})=O(n^{1+1/\ell})$. This upper bound is only known to be tight for $\ell=2,3$ and $5$. For more results on bipartite-graph Tur\'an problems, see the informative survey by F\"uredi and Simonovits~\cite{furedi2013history}. For hypergraphs, less is known. Tur\'an~\cite{turan1941} conjectured that $\ex(n,K^3_4)=(\frac{5}{9}+o(1))\binom{n}{3}$ and $\ex(n,K^3_5)=(\frac{3}{4}+o(1))\binom{n}{3}$, where $K^3_t$ is the complete $3$-graph on $t$ vertices. This conjecture, though innocent-looking, has been open for more than 80 years. As for the subject of this paper, linear cycles, F\"uredi and Jiang~\cite{furedi2014hypergraph} and Kostochka, Mubayi and Verstra\"ete~\cite{kostochka2015turan} showed that $\ex(n,C^r_{\ell})=(\lfloor\frac{\ell-1}{2}\rfloor+o(1))\binom{n}{r-1}$. See~\cite{keevash2011hypergraph} for an excellent survey by Keevash on hypergraph Tur\'an problems.

Let $G_{n,p}^r$ denote the random $r$-graph where edges of $K_n^{r}$ are sampled independently with probability $p$. The random Tu\'ran number $\ex(G^r_{n,p},H)$ is the maximum number of edges in an $H$-free subgraph of $G^r_{n,p}$. The random Tu\'ran problem consists of estimating the random variable $\ex(G^r_{n,p},H)$. To this end, we say that a statement depending on $n$ holds {\em asymptotically almost surely} (abbreviated a.a.s.) if the probability that it holds tends to 1 as $n$ tends to infinity. See~\cite{rodl2013extremal} for a great survey by R\"odl and Schacht on random Tur\'an problems. Since $\ex(G^r_{n,1},H)=\ex(n,H)$, it is natural to ask whether $\ex(G^r_{n,p},H)=(1+o(1))p\cdot\ex(n,H)$ a.a.s. when $p$ is large. This is true for $H$ that is not $r$-partite due to the breakthrough work done independently by Conlon and Gowers~\cite{conlon2016combinatorial} and Schacht~\cite{schacht2016extremal}. To formally state their result, we define the {\em $r$-density} of $H$,
$$
m_r(H)=\max_{G\subset H,~e(G)\ge 2}\l\{\frac{e(G)-1}{v(G)-r}\r\}.
$$
If $p\ll n^{-1/m_r(H)}$, then the number of copies of $H$ in $G^r_{n,p}$ is much smaller than the number of edges in $G^r_{n,p}$, hence $\ex(G^r_{n,p},H)=(1+o(1)e(G^r_{n,p})$. Conlon and Gowers~\cite{conlon2016combinatorial} and Schacht~\cite{schacht2016extremal} showed the following, which confirm conjectures of Kohayakawa, {\L}uczak and R\"odl~\cite{kohayakawa1997k} and R\"odl, Ruci\'nski and Schacht~\cite{rodl2007ramsey}.
\begin{thm}[\cite{conlon2016combinatorial,schacht2016extremal}]
     If $p\gg n^{-1/m_r(H)}$, then as $n\to\infty$, a.a.s.
     $$
     \ex(G^r_{n,p},H)=p(\ex(n,H)+o(n^r)).
     $$
\end{thm}

The behaviour of $\ex(G^r_{n,p},H)$ when $H$ is $r$-partite is still a wide open problem. There are many sporadic results, see for examples~\cite{furedi1994random,haxell1995turan,MORRIS2016534, balogh2011number,jiang2022balanced,spiro2022random,mubayi2020random,nie2021triangle,spiro2021relative,spiro2022counting}. One extensively studied case is when $H=C_{2\ell}$, the graph cycle of length $2\ell$. Haxell, Kohayakawa and {\L}uczak~\cite{haxell1995turan} showed that a.a.s. $\ex(G^2_{n,p}, C_{2\ell})\ll e(G^2_{n,p})$ when $p\gg n^{-1/m_2(C_{2\ell})}=n^{-1+\frac{1}{2\ell-1}}$. Using the container method, Morris and Saxton~\cite{MORRIS2016534} further improved the upper bounds.

\begin{thm}[\cite{haxell1995turan, MORRIS2016534}]\label{theorem:morris-saxton}
    For every $\ell\ge 2$, there exists $C=C(\ell)$ such that a.a.s.
    $$
    \ex(G^2_{n,p},C_{2\ell})\le
    \l\{
    \begin{aligned}
    &n^{1+\frac{1}{2\ell-1}}(\log n)^2,~&\text{if}~n^{-1+\frac{1}{2\ell-1}}\ll p\le n^{-\frac{\ell-1}{2\ell-1}}(\log n)^{2\ell};\\    
    &Cp^{\frac{1}{\ell}}n^{1+\frac{1}{\ell}},&\text{if}~p\ge n^{-\frac{\ell-1}{2\ell-1}}(\log n)^{2\ell}.
    \end{aligned}
    \r.
    $$
\end{thm}

This theorem is essentially tight conditional on the following notorious conjecture of Erd\H{o}s and Simonovits~\cite{erdHos1982compactness}.
\begin{conj}[\cite{erdHos1982compactness}]\label{conjecture:Erdos-Simonovits}
For $\ell\ge2$,
    $$\ex(n,\{C_3,\dots,C_{2\ell}\})=\Theta(n^{1+1/\ell}).$$
\end{conj}

The $r$-uniform linear cycle of length $\ell\ge 3$, denoted by $C^r_\ell$, is an $r$-uniform hypergraph with $\ell$ edges $e_1,\dots,e_\ell$ and $\ell(r-1)$ vertices $v_1,\dots,v_{\ell(r-1)}$ such that $e_i=\{v_{(i-1)(r-1)},\dots,v_{i(r-1)}\}$, $\forall 1\le i\le \ell$ (here $v_0=v_{\ell(r-1)}$). As an attempt to extend Theorem~\ref{theorem:morris-saxton} to hypergraphs, Mubayi and Yepremyan~\cite{mubayi2020random} showed the following.
\begin{thm}[\cite{mubayi2020random}]\label{theorem:mubayi-yepremyan}
    For $r\ge 3$ and $\ell\ge 2$ a.a.s. the following holds:
    $$
    \ex(G^r_{n,p},C^r_{2\ell})\le
    \l\{
    \begin{aligned}
        &p^{\frac{1}{2\ell-1}}n^{1+\frac{r-1}{2\ell-1}+o(1)},~~&\text{if}~n^{-r+2+o(1)}\le p\le n^{-r+2+\frac{1}{2\ell-2}+o(1)};\\
        &pn^{r-1+o(1)},~~&\text{if}~ p \ge n^{-r+2+\frac{1}{2\ell-2}+o(1)}.
    \end{aligned}
    \r.
    $$
\end{thm}

This upper bound is not tight when $n^{-r+2+o(1)}\le p\le n^{-r+2+\frac{1}{2\ell-2}+o(1)}$. In this range, it is known~\cite{mubayi2020random} that 
\begin{equation}\label{equation:lower}
    \ex(G^r_{n,p},C^r_{2\ell})\ge \max\{n^{1+\frac{1}{2\ell-1}+o(1)},\Omega(pn^{r-1})\}.
\end{equation}

The authors of~\cite{mubayi2020random} conjectured that (\ref{equation:lower}) is tight for $C^3_4$.

\begin{conj}[\cite{mubayi2020random}]
$$
\ex(G^3_{n,p},C^3_4)=
\l\{
\begin{aligned}
    &(1+o(1))e(G^3_{n,p}),~~~~~&\text{if}~n^{-3}\ll p\ll n^{-5/3};\\
    &\Theta(n^{4/3+o(1)}),~~~~~&\text{if}~n^{-5/3}\ll p\ll n^{-2/3};\\
    &\Theta(pn^2),~~~~~&\text{otherwise}.
\end{aligned}
\r.
$$
\end{conj}

Theorem~\ref{theorem:mubayi-yepremyan} relies on a new balanced supersaturation result for $C^r_{2\ell}$. ``Balanced supersaturation'' roughly means that, given an $r$-graph $H$, if the number of edges in an $r$-graph $G$ on $n$ vertices is much larger than $\ex(n,H)$, then we can find a collection $\S$ of copies of $H$ in $G$ with ``balanced distribution'', i.e. one cannot find a set $\sigma$ of edges of $G$ such that the number of copies of $H$ in $\S$ containing $\sigma$ is unusually large. A good balanced supersaturation result for $H$ would usually imply a good upper bound for the random Tur\'an number of $H$ by using the container method developed independently by Balogh, Morris and Samotij~\cite{balogh2015independent} and Saxton and Thomason~\cite{saxton2015hypergraph}. 

Let $\Fb(n,H)$ be the number of $H$-free hypergraphs on $n$ vertices. Determining $\Fb(n,H)$ is another type of problem that relies heavily on balanced supersaturation results. This type of problem has been studied extensively for $r$-graphs that are not $r$-partite~\cite{erdos1986asymptotic,nagle2001asymptotic,nagle2006extremal}. Mubayi and Wang~\cite{mubayi2019thenumber} initiated the study for $r$-partite $r$-graphs. In~\cite{mubayi2019thenumber}, they determined the asymptotics of $\Fb(n,C^r_k)$ for even $k$ and $r=3$, and conjectured the asymptotics for all $r$ and $k$. Their conjecture was later confirmed by Balogh, Skokan and Narayanan~\cite{balogh2019number}. Ferber, Mckinley and Samotij~\cite{ferber2020supersaturated} proved similar results for a larger family of hypergraphs which includes linear cycles. The results of~\cite{balogh2019number} and \cite{ferber2020supersaturated} both rely on balanced supersaturation. More recently, Jiang and Longbrake~\cite{jiang2022balanced} proved the balanced supersaturation result in~\cite{ferber2020supersaturated} in a more explicit form.  However, the balanced supersaturation results mentioned above are not strong enough to imply tight upper bounds for the random Tur\'an numbers of linear cycles.

\subsection{Our results}

In this paper, we show that (\ref{equation:lower}) is tight for all $r\ge 4$ using ideas inspired by Mubayi and Yepremyan~\cite{mubayi2020random} and Balogh, Skokan and Narayanan~\cite{balogh2019number}. The two papers used very different strategies to obtain balanced supersaturation results: the strategy of \cite{mubayi2020random} works better when the ``typical'' codegree is ``small''; and the strategy of \cite{balogh2019number} outperforms when the ``typical'' codegree is ``large''. We manage to combine these two ideas using codegree dichotomy.
\begin{thm}\label{theorem:main_ge4}
    For every $r\ge 4$ and $\ell\ge2$, there exists $c>0$ such that if $p\ge (\log n)^{2r-2}n^{-r+2+\frac{1}{2\ell-1}}$, then a.a.s.
    $$
    \ex(G^r_{n,p}, C^r_{2\ell})\le cpn^{r-1}.
    $$
\end{thm}

When $p\le n^{-r+2+\frac{1}{2l-1}+o(1)}$, since $\ex(G^r_{n,p}, C^r_{2\ell})$ can be (more or less) considered as a non-decreasing function with respect to $p$, we have a.a.s.
$$
    \ex(G^r_{n,p}, C^r_{2\ell})\le n^{1+\frac{1}{2\ell-1}+o(1)}.
$$
Due to (\ref{equation:lower}), this is tight when $n^{-r+2+o(1)}\le p\le n^{-r+2+\frac{1}{2l-1}+o(1)}$.

Note that $p=n^{-r+2+\frac{1}{2l-1}+o(1)}$ is the point where $\Omega\l(pn^{r-1}\r)$ becomes the maximum in (\ref{equation:lower}). Hence Theorem~\ref{theorem:main_ge4}, together with (\ref{equation:lower}), close the gap for $\ex(G^r_{n,p}, C^r_{2\ell})$ in the range when $n^{-r+2+o(1)}\le p\le n^{-r+2+\frac{1}{2\ell-2}+o(1)}$.

The author recently learned that Mubayi and Yepremyan~\cite{Mubayi2023OnTR} have also obtained essentially the same result as Theorem~\ref{theorem:main_ge4} in a different way. They followed the same strategy as~\cite{mubayi2020random} with more careful estimation. We will briefly discuss the difference of our methods in Subsection~\ref{section:outline}.

\begin{figure}[h]
    \centering
    \includegraphics[scale=0.4]{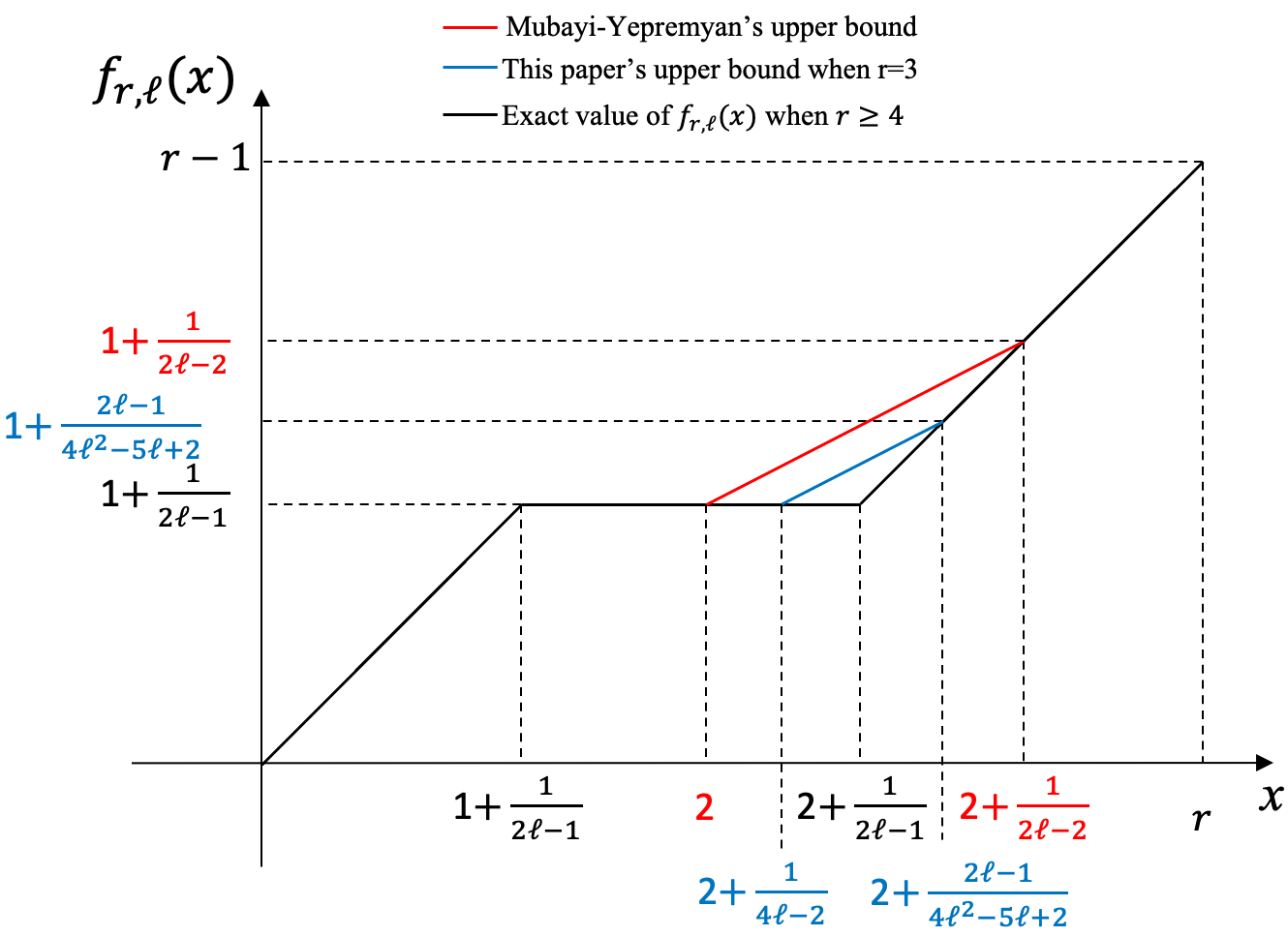}
    \caption{The behaviour of $f_{r,\ell}(x)$ with respect to $x$}
    \label{fig:plot}
\end{figure}

For $r=3$, we substantially improve Theorem~\ref{theorem:mubayi-yepremyan} (See Figure~\ref{fig:plot}), but unfortunately it is still not tight. 
\begin{thm}\label{theorem:main_3}
    For every $\ell\ge2$, there exists $c>0$ such that the following holds. Let
    $$
    p_0=cn^{-\frac{4\ell-3}{4\ell-2}}(\log n)^{\frac{4\ell^2+\ell-4}{2\ell-2}},~p_1=n^{-\frac{(\ell-1)(4\ell-3)}{4\ell^2-5\ell+2}}(\log n)^{3+\frac{2\ell-2}{4\ell^2-5\ell+2}}.
    $$
    Then a.a.s.
    $$
    \ex(G^3_{n,p}, C^3_{2\ell})\le 
    \l\{
    \begin{aligned}
        &cp^{\frac{2(\ell-1)}{\ell(4\ell-3)}}n^{1+\frac{1}{\ell}}(\log n)^{3-\frac{4\ell-4}{\ell(4\ell-3)}},~~&\text{if}~p_0\le p<p_1;\\
        &cpn^2,~~&\text{if}~p\ge p_1.
    \end{aligned}
    \r.
    $$
\end{thm}

Together with simple lower bounds, the theorems above imply the following corollaries, which describe roughly the behavior of $\ex(G^r_{n,p}, C^r_{2\ell})$ with respect to $p$ in a simpler way. The limits of random variables below are defined with convergence in distribution.

\begin{cor}
For integers $r\ge4$, $\ell\ge 2$ and real number $0< x\le r$, let $p=n^{-r+x}$ and define
\begin{equation*}
f_{r,\ell}(x)=\lim_{n\rightarrow\infty}\log_n\ex(G^{r}_{n,p}, C^{r}_{2\ell}).
\end{equation*}
Then we have
\begin{equation*}
f_{r,\ell}(x)=\l\{
    \begin{aligned}
    &x,~~~~~~~~~&0\le x\le 1+\frac{1}{2\ell-1},\\
    &1+\frac{1}{2\ell-1},~~~&1+\frac{1}{2\ell-1}\le x\le 2+\frac{1}{2\ell-1},\\
    &x-1,~~~~~~~~~ &2+\frac{1}{2\ell-1}\le x\le r.\\
    \end{aligned}
    \r.
\end{equation*}
\end{cor}

For $r=3$, the same lower bound holds, but the upper bound is worse for a small range of $x$.

\begin{cor}
For integer $\ell\ge 2$ and real number $0< x\le 3$, let $p=n^{-3+x}$ and define
\begin{equation*}
f_{3,\ell}(x)=\lim_{n\rightarrow\infty}\log_n\ex(G^{3}_{n,p}, C^{3}_{2\ell}).
\end{equation*}
Then we have
\begin{equation*}
f_{3,\ell}(x)=\l\{
    \begin{aligned}
    &x,~~~~~~~~~&0\le x\le 1+\frac{1}{2\ell-1},\\
    &1+\frac{1}{2\ell-1},~~~&1+\frac{1}{2\ell-1}\le x\le 2+\frac{1}{4\ell-2},\\
    &x-1,~~~~~~~~~ &2+\frac{2\ell-1}{4\ell^2-5\ell+2}\le x\le 3.\\
    \end{aligned}
    \r.
\end{equation*}
When $2+\frac{1}{4\ell-2}\le x\le 2+\frac{2\ell-1}{4\ell^2-5\ell+2}$, 
$$
\max\{1+\frac{1}{2\ell-1}, x-1\}\le f_{3,\ell}(x)\le \frac{2(\ell-1)}{\ell(4\ell-3)}x+\frac{4\ell^2-5\ell+3}{\ell(4\ell-3)}.
$$
\end{cor}

Here we roughly describe the constructions for the lower bounds. 

When $x<1+\frac{1}{2\ell-1}$, the number of copies of $C^r_{2\ell}$ is much smaller than the number of edges in the random hypergraph, hence with high probability we can obtain a $C^r_{2\ell}$-free subgraph with $pn^{r-o(1)}=n^{x-o(1)}$ edges by deleting an edge in each copy of $C^r_{2\ell}$. 

When $1+\frac{1}{2\ell-1}\le x\le 2+\frac{1}{2\ell-1}$, we pick a subgraph $H$ of $G^r_{n,p}$ such that edges are sampled independently with probability $p'=n^{-r+1+\frac{1}{2\ell-1}-o(1)}/p$. Clearly, $H=G^r_{n,pp'}$. It is not hard to check that with high probability $e(H)$ is at least $n^{1+\frac{1}{2\ell-1}-o(1)}$ and the number of copies of $C^r_{2\ell}$ in $H$ is much smaller than $e(H)$. Hence we can obtain a $C^r_{2\ell}$-free subgraph with $n^{1+\frac{1}{2\ell-1}-o(1)}$ edges by deleting an edge in each copy of $C^r_{2\ell}$.

When $x\ge 2+\frac{1}{2\ell-1}$, we fix a vertex $v$ and pick a subgraph $H$ of $G^r_{n,p}$ whose edge set consists of all edges of $G^r_{n,p}$ containing $v$. Clearly, $H$ is $C^r_{2\ell}$-free. Moreover, with high probability, $e(H)\ge pn^{r-1-o(1)}=n^{x-1-o(1)}$.

%BERGE

Our method can also be used to give upper bounds for the random Tur\'an number of Berge even cycles. For $r\ge k\ge 2$, the $r$-uniform {\em Berge cycles} of length $k$, denoted by $\B^r_k$, is the family of $r$-graphs with $k$ edges $e_1,\dots,e_k$ such that there exist $k$ vertices $v_1,\dots,v_k$ with the property that $v_{i-1}, v_i\in e_i$ for all $1\le i\le k$, where $v_0=v_k$. Further, we let $\B^r_{[k]}=\cup_{2\le i\le k}\B^r_i$. Spiro and Verstra\"ete~\cite{spiro2022counting} showed the following. 

\begin{thm}[\cite{spiro2022counting}]\label{theorem:Spiro-Vers}
    Let $\ell \ge 2$ and $r\ge 2$, and let $\lambda=\lceil(r-2)/(2\ell-2)\rceil$. Then a.a.s.:
    $$
    \ex(G^r_{n,p},\B^r_{[2\ell]})\le
    \l\{
    \begin{aligned}
        &n^{1+\frac{1}{2\ell-1}+o(1)},~&\text{if}~n^{-r+1+\frac{1}{2\ell-1}}\le p\le n^{-\frac{(r-1+\lambda)(\ell-1)}{2\ell-1}+o(1)};\\
        &p^{\frac{1}{(r-1+\lambda)\ell}}n^{1+\frac{1}{\ell}+o(1)},~&\text{if}~n^{-\frac{(r-1+\lambda)(\ell-1)}{2\ell-1}+o(1)}\le p\le 1.
    \end{aligned}
    \r.
    $$
    If Conjecture~\ref{conjecture:Erdos-Simonovits} is true, then
    $$
    \ex(G^r_{n,p},\B^r_{[2\ell]})\ge
    \l\{
    \begin{aligned}
        &n^{1+\frac{1}{2\ell-1}+o(1)},~&\text{if}~n^{-r+1+\frac{1}{2\ell-1}}\le p\le n^{-\frac{(r-1)(\ell-1)}{2\ell-1}};\\
        &p^{\frac{1}{(r-1)\ell}}n^{1+\frac{1}{\ell}-o(1)},~&\text{if}~n^{-\frac{(r-1)(\ell-1)}{2\ell-1}}\le p\le 1.
    \end{aligned}
    \r.
    $$
\end{thm}
We note that the proof of Theorem~\ref{theorem:main_ge4} and Theorem~\ref{theorem:main_3} can be easily adapted to show the following theorem which improves Theorem~\ref{theorem:Spiro-Vers}.
\begin{thm}\label{Theorem:Berge}
Let $\ell \ge 2$ and $r\ge 2$, and let $\lambda=(r-2)/(2\ell-2)$. Then a.a.s.:
    $$
    \ex(G^r_{n,p},\B^r_{2\ell})\le
    \l\{
    \begin{aligned}
        &n^{1+\frac{1}{2\ell-1}+o(1)},~&\text{if}~n^{-r+1+\frac{1}{2\ell-1}}\le p\le n^{-\frac{(r-1+\lambda)(\ell-1)}{2\ell-1}+o(1)};\\
        &p^{\frac{1}{(r-1+\lambda)\ell}}n^{1+\frac{1}{\ell}+o(1)},~&\text{if}~n^{-\frac{(r-1+\lambda)(\ell-1)}{2\ell-1}+o(1)}\le p\le 1.
    \end{aligned}
    \r.
    $$
\end{thm}
Theorem~\ref{Theorem:Berge} improves Theorem~\ref{theorem:Spiro-Vers} in the following sense: the $\lambda$ in Theorem~\ref{Theorem:Berge} does not have the ceiling operation; and the $\B^{r}_{[2\ell]}$ in Theorem~\ref{theorem:Spiro-Vers} is replaced by $\B^r_{2\ell}$ in Theorem~\ref{Theorem:Berge}. However, there is still a gap between our upper bounds and lower bounds in Theorem~\ref{theorem:Spiro-Vers}. It is not clear to the author whether the upper bounds or lower bounds are closer to the truth.

\subsection{Proof outline}\label{section:outline}
The results in this paper mainly come from proving stronger balanced supersaturation theorems for $C^r_{2\ell}$. Obtaining upper bounds on random Tur\'an numbers from balanced supersaturation results is by far a routine application of the container method developed independently by Balogh, Morris and Samotij~\cite{balogh2015independent} and Saxton and Thomason~\cite{saxton2015hypergraph}. Hence we will only outline the proof of the new balanced supersaturation theorems. 

Given an $r$-graph $H$, recall that ``balanced supersaturation'' roughly means that if the number of edges in an $r$-graph $G$ on $n$ vertices is much larger than $\ex(n,H)$, then we can find a collection $\S$ of copies of $H$ in $G$ with ``balanced distribution''. To construct a balanced collection of copies of $C^r_{2\ell}$, we will adopt the following strategy using codegree dichotomy and induction on $r$. 

For the base case $r=2$, we make use of the balanced supersaturation theorems for $C_{2\ell}$ by Morris and Saxton~\cite{MORRIS2016534}. Now suppose we already have a good balanced supersaturation theorem for $C^{r-1}_{2\ell}$. Given an $r$-graph $H$ with many edges, we can assume that it is $r$-partite with respect to the partition $V(H)=V_1\cup\cdots\cup V_r$. A {\em $k$-shadow} of $H$ is a set $\sigma$ of $k$ vertices such that $\sigma\subset e$ for some edge $e$ of $H$. The {\em codegree} of $\sigma$ is the number of edges of $H$ containing $\sigma$. The proof splits into two cases according to how large is the ``typical'' codegree of the $(r-1)$-shadows of $H$. More precisely, $H$ can be partitioned into subgraphs $F$ and $F_{\tau,a}$, $\tau\in[r]^{r-1}$, $a\ge 0$, such that: 1. the codegree of every $(r-1)$-shadow of $F$ is ``large'', and 2. every $(r-1)$-shadow in $\cup_{i\in\tau}V_i$ of $F_{\tau,a}$ has codegree at least $2^a$ and at most $2^{a+1}$ where $2^a$ is ``small''. The idea of partitioning $H$ into subgraphs with ``almost regular'' codegrees has already appeared in~\cite{mubayi2020random}. What we do here is slightly different: we partition the $r$-graph according to the codegree of its $(r-1)$-shadows, while they do it according to the codegree of its $2$-shadows; moreover, we only require that the codegrees of $(r-1)$-shadows in $\cup_{i\in\tau}V_i$ are ``almost regular'' for SOME $\tau\in[r]^{r-1}$, while they require that the codegrees of $2$-shadows in $V_i\cup V_j$ are ``almost regular'' for ALL pairs $(i,j)$.

Now if $H$ has a subgraph $F$ with $e(H)/(\log n)^{r-1}$ edges such that all $(r-1)$-shadows of $F$ have ``large'' codegree, then we can construct a balanced collection of $C^r_{2\ell}$ using ``greedy expansion''. For example, to greedily construct a copy of $C^3_4$ with edges $\{v_1,w_4,w_1\}$, $\{v_2,w_1,w_2\}$, $\{v_3,w_2,w_3\}$ and $\{v_4,w_3,w_4\}$, we first arbitrarily pick an edge $e_1=\{v_1,w_4,w_1\}$, and then greedily pick edges $e_2=\{w_4,w_1,w_3\}$, $e_3=\{w_1,w_3,w_2\}$, $e_4=\{v_2,w_1,w_2\}$, $e_5=\{v_3,w_2,w_3\}$ and $e_6=\{v_4,w_3,w_4\}$. Because the codegree of every $2$-shadow is ``large'', we can obtain many $C^3_4$ in this way; with some additional techniques we can also ensure that they are evenly distributed. This ``greedy expansion'' method is inspired by Balogh, Narayanan and Skokan~\cite{balogh2019number}. 

On the other hand, if such $F$ does not exist, then we can find $(r-1)$ parts, say $V_1,\dots,V_{r-1}$, and a subgraph $F'$ of $H$ with $e(F')\ge e(H)/\log n$ such that all $(r-1)$-shadows of $F'$ in $V_1\cup\cdots\cup V_{r-1}$ have codegree at least $D$ and at most $2D$ where $D$ is ``small''. We now consider the $(r-1)$-graph $G$ on $V_1\cup\dots\cup V_{r-1}$ whose edge set consists of all $(r-1)$-shadows of $F'$ in $V_1\cup\cdots\cup V_{r-1}$. Then $e(G)\ge e(F')/(2D)$. Since $D$ is ``small'', $e(G)$ must be ``large''. Applying the balanced supersaturation theorem for $C^{r-1}_{2\ell}$ on $G$ gives a balanced collection $\S'$ of copies of $C^{r-1}_{2\ell}$ in $G$. Using the fact that the codegrees are ``almost regular'', we can extend each copy of $C^{r-1}_{2\ell}$ in $\S'$ into copies of $C^r_{2\ell}$ in $H$, which gives a balanced collection $\S$ of copies of $C^r_{2\ell}$ in $H$. This idea of first using the supersaturation in the shadow and then extending can be traced back to Mubayi and Yepremyan~\cite{mubayi2020random}.

By selecting the best cutoff that defines ``large'' and ``small'' codegrees, we obtain the supersaturation theorems for $C^r_{2\ell}$ in this paper. 

To summarize, given an $r$-graph $H$ with many edges, if the typical codegrees of its $(r-1)$-shadows are ``large'', then we find copies of $C^r_{2\ell}$ by ``greedy expansion''; if the typical codegrees are ``small'', then we find copies of $C^{r-1}_{2\ell}$ in its $(r-1)$-shadows and then expand them into copies of $C^r_{2\ell}$. The previous result of Mubayi and Yepremyan~\cite{mubayi2020random} adopted a similar but different approach: Given an $r$-graph $H$ with many edges, they find many $2$-shadows of $H$ with ``small'' and ``almost regular'' codegrees. Then they make use of the supersaturation theorem for $C_{2\ell}$ to find copies of $C_{2\ell}$ in the aforementioned set of $2$-shadows of $H$, and then expand them into copies of $C^r_{2\ell}$ in $H$. The reason why we have a better result for $r=3$ is that the ``small'' codegree in our proof is smaller than the ``small'' codgree in their proof. For $r\ge 4$, Mubayi and Yepremyan~\cite{Mubayi2023OnTR} manage to show that the ``typical'' codegrees in some $V_i\cup V_j$ are ``smaller'' than that in~\cite{mubayi2020random}, which gives them the improvement. In this case, our method and theirs perform equally well.

\subsection{Notations and structure}
Given an $r$-graph $H=(V,E)$, we denote $|V|$ by $v(H)$ and $|E|$ by $e(H)$. A {\em $k$-shadow} of $H$ is a set $\sigma$ of $k$ vertices such that $\sigma\subset e$ for some edge $e$ of $H$. Given a set $\sigma$ of vertices of $H$, let $N_H(\sigma)$ be the set of edges of $H$ containing $\sigma$ and let $d_{H}(\sigma)=|N_H(\sigma)|$. For $1\le j\le r$, let $\Delta_j(H)$ denote the maximum $d_H(\sigma)$ with $|\sigma|=j$. Given a collection $\S$ of copies of $C^r_{2\ell}$ in an $r$-graph $H$, we can view $\S$ as a hypergraph on $E(H)$ whose edges are copies of $C^r_{2\ell}$ in $\S$. Then $\Delta_j(\S)$ means the maximum number of copies of $C^r_{2\ell}$ in $\S$ containing $j$ common edges. Given two positive functions $f(n)$ and $g(n)$ on $\N$, we write $f(n)\le O(g(n))$ if there exists a constant $c>0$ such that $f(n)\le c\cdot g(n)$ for all $n$, and write $f(n)= \Omega(g(n))$ if $g\le O(f)$. If $f= O(g)$ and $g= O(f)$, then we write $f=\Theta(g)$. We write $f(n)=o(g(n))$ if $f(n)/g(n)\to 0$ as $n\to\infty$, and write $f(n)=\omega(g(n))$ if $g(n)=o(f(n))$.  

The rest of this paper will be organized as follows: In Section~\ref{section:supersaturation}, we prove the balanced supersaturation theorems for linear cycles of even length. In Section~\ref{section:container}, we make use of the container method together with the supersaturation theorems to obtain container theorems for linear cycles of even length. More precisely, we find a ``small'' collection of $r$-graphs on $[n]$ with ``few'' edges such that they contain all $C^r_{2\ell}$-free $r$-graphs on $[n]$. In Section~\ref{section:upper}, we prove Theorem~\ref{theorem:main_ge4} and Theorem~\ref{theorem:main_3}, the random Tur\'an theorems for linear even cycles. In Section~\ref{section:berge}, we briefly discuss why the proof of Theorem~\ref{theorem:main_ge4} and Theorem~\ref{theorem:main_3} can be adapted into a proof of Theorem~\ref{Theorem:Berge}, the random Tur\'an theorem for Berge even cycles.
%%%%%%%%%%%%%%%%%%%%%%%%%%%%%%%%%%%%%%%%%%%%%%%%%%%%%%%%%%%%%%%%%%%%%%%%%%%%%%%%%%%%%%%%%%%%%%%%%%%%%%%%%%%%%%%%%%%%%%%%%%%%%%%%%%%%%%%%%%%%%%%%%%%%%%%%%%%%%%%%%%%%%%%%%%%%%%%%%%%%%%%%%%%%%%%%%%%%%%%%%%%%%%%%%%%%%%%%%%%%%%%%%%%%%%%%%%%%%%%%%%%%%%%%%%%%%%%%%%%%%%%%%%%%%%%%%%%%%%%%%%%%%%%%%%%%%%%%%%%%%%%%%%%%%%%%%%%%%%%%

\section{Balanced supersaturation}\label{section:supersaturation}

In this section, we prove the balanced supersaturation theorems for $C^r_{2\ell}$. We use the following balanced supersaturation theorem for $C_{2\ell}$ by Morris and Saxton~\cite{MORRIS2016534} as the base case of our inductive proof.

\begin{thm}[Morris and Saxton~\cite{MORRIS2016534}]\label{theorem:cycle_supersaturation}
    For every $\ell\ge 2$, there exist $c, K>0$ such that the following holds for every $t\ge K$ and every $n$ sufficiently large. Given a graph $G$ with $n$ vertices and $tn^{1+1/\ell}$ edges there exists a collection $\S$ of copies of $C_{2\ell}$ in $G$ satisfying
%     \begin{enumerate}
%         \item[(a)]$|\S|\ge ct^{2\ell}n^2$,
%         \item[(b)]$\Delta_{k}(\S)\le c^{-1}t^{2\ell-k-\frac{k-1}{l-1}}n^{1-\frac{1}{l}}$,~$\forall 1\le k\le 2\ell-1$.
%     \end{enumerate}
% As a result of the two inequalities above, we have
$$
\Delta_j(\S)\le\frac{c|\S|}{tn^{1+1/\ell}}\l(\max\l\{t^{-\frac{\ell}{\ell-1}},~t^{-1}n^{-\frac{\ell-1}{\ell(2\ell-1)}}\r\}\r)^{j-1},~\forall 1\le j\le 2\ell.
$$
\end{thm}

Theorem~\ref{theorem:cycle_supersaturation} is slightly weaker than the original Theorem 1.5 in~\cite{MORRIS2016534} -- Here we state the balanced supersaturation result more compactly, which makes it easier for the later application of the container theorem. Note that the term $t^{-1}n^{-\frac{\ell-1}{\ell(2\ell-1)}}$ comes from the fact that $\Delta_{2\ell}(\S)=1$.

The following result is folklore.

\begin{prop}[Balogh, Narayanan and Skokan~\cite{balogh2019number}]\label{prop:regular}
    For any integers $n>t>0$, there exists a graph on $n$ vertices such that each vertex has degree $t$ or $t-1$.
\end{prop}

The following lemma is inspired by Theorem 3.1 in~\cite{balogh2019number} of Balogh, Narayanan and Skokan. 

\begin{lem}\label{lem: expansion}
    For every $r\ge 3$ and $\ell\ge 2$, there exist constants $K, c>0$ such that the following holds for all sufficiently large integer $n$ and every $t\ge K$. Let $H$ be an $r$-graph on $n$ vertices. If every $(r-1)$-shadow $\sigma$ of $H$ has $d_H(\sigma)\ge t$, then there exists a collection $\S$ of copies of $C^r_{2\ell}$ in $H$ satisfying
    \begin{enumerate}
        \item[(a)] $|\S|\ge ce(H)t^{(2\ell-1)r-2\ell}$;
        \item[(b)] $\Delta_j(\S)\le c^{-1}t^{(2\ell-j)(r-1)-1},~1\le j\le 2\ell-1$.
    \end{enumerate}
\end{lem}

\begin{proof}
    For each $(r-1)$-shadow $\sigma$ in $H$, since $d_H(\sigma)\ge t$, by Proposition~\ref{prop:regular} we can define a graph $\Gamma(\sigma)$ on $N_H(\sigma)$ such that, for each $e\in N_H(\sigma)$, $d_{\Gamma(\sigma)}(e)=t$ or $t-1$. Suppose $K\ge 4\ell(r-1)$. We now describe an algorithm to construct copies of $C^r_{2\ell}$ in $H$; it involves specifying $2\ell$ edges $e_1,\dots,~e_{2\ell}$ which form a copy of $C^r_{2\ell}$ in $H$, where $e_i=\{v_{i,1},\dots,~v_{i,r-2},~w_{i-1},~w_{i}\}$, $1\le i\le 2\ell$, $w_0=w_{2\ell}$. In this algorithm, we can view an edge as an ordered set of vertices. Given an edge with an order on its vertices $e=(a_1,~a_2,\dots,~a_r)$ and a vertex $v$, we let $e[2,r]=(a_2,~\dots,~a_r)$ and let $e[2,r]\uplus\{v\}=(a_2,~\dots,~a_r,~v)$.
    \begin{enumerate}
        \item[(i)] We start by choosing an edge with an order $e_1=(v_{1,1},\dots,~v_{1,r-2},w_{2\ell},~w_1)$ in $E(H)$; 
        \item[(ii)] For $1\le i\le\ell-1$, we specify $w_{2\ell-i}$ and $w_{i+1}$ inductively as following. Let $f_1=e_1$. Next we greedily pick $f_{2\ell}=f_1[2,r]\uplus\{w_{2\ell-1}\}$ such that $f_{2\ell}$ and $f_1$ are adjacent in $\Gamma(f_1[2,r])$, and that $w_{2\ell-1}$ is different from any specified vertices. Then we pick $f_{2}=f_{2\ell}[2,r]\uplus\{w_2\}$ such that $f_{2}$ and $f_{2\ell}$ are adjacent in $\Gamma(f_1[2,r])$ and that $w_2$ is different from any other specified vertices. In general, given $f_i$, we choose a neighbor $f_{2\ell+1-i}=f_i[2,r]\uplus\{w_{2\ell-1}\}$ of $f_i$ in $\Gamma(f_i[2,r])$ such that the only vertex $w_{2\ell-i}$ in $f_{2\ell+1-i}\setminus f_i$ is distinct from all specified vertices. Then we choose a neighbor $f_{i+1}=f_{2\ell+1-i}[2,r]\uplus\{w_{i+1}\}$ of $f_{2\ell+1-i}$ in $\Gamma(f_{2\ell+1-i}[2,r])$ such that the only vertex $w_{i+1}$ in $f_{i+1}\setminus f_{2\ell+1-i}$ is distinct from all specified vertices. Let $f_{\ell+1}=f_\ell$. So far we have $\{w_{i-1},~w_i\}\subset f_i$ for every $2\le i\le 2\ell$. (See Figure~\ref{fig:greedy expansion} for an example.)
    \begin{figure}[h]
        \centering
        \includegraphics[scale=0.3]{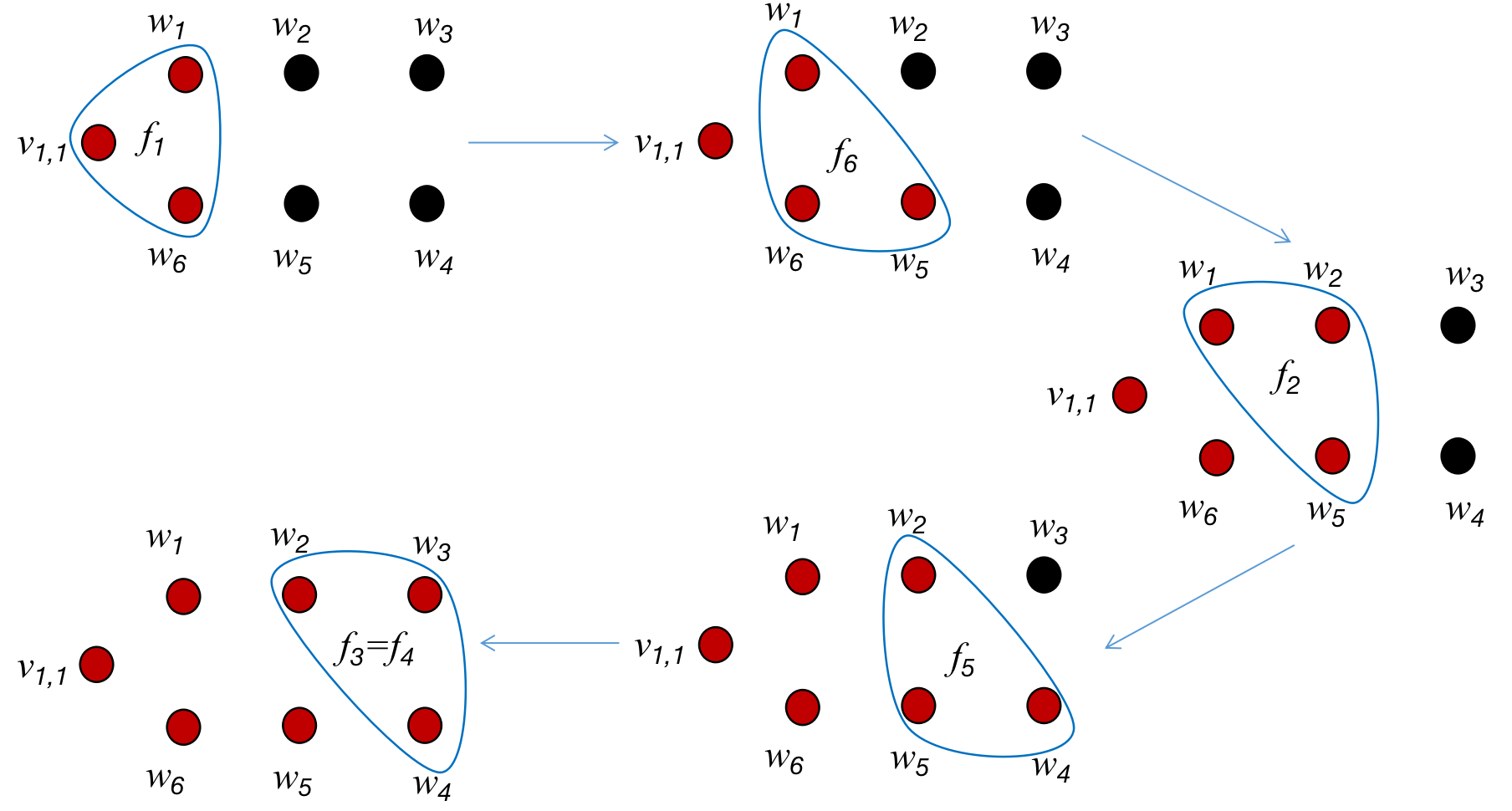}
        \caption{Step (ii) of the algorithm for $C^3_6$: red vertices are specified.}
        \label{fig:greedy expansion}
    \end{figure}
        \item[(iii)] For $2\le i\le 2\ell$, we give $f_i$ a new order such that $w_{i-1}$ and $w_i$ are the last two vertices. Then we specify $v_{i,j}$, $1\le j\le r-2$, inductively as following. Let $e_{i,0}=f_i$. We choose a neighbor $e_{i,j}$ of $e_{i,j-1}$ in $\Gamma(e_{i,j-1}[2,r])$ such that the only vertex $v_{i,j}$ in $e_{i,j}\setminus e_{i,j-1}$ is distinct from all specified vertices; give vertices of $e_{i,j}$ an order such that $e_{i,j}=e_{i,j-1}[2,r]\uplus\{v_{i,j}\}$. Let $e_{i}=e_{i,r-2}$.
    \end{enumerate}
    It is not hard to check that the $2\ell$ edges $e_1,\dots,e_{2\ell}$ generated by the algorithm form a copy of $C^r_{2\ell}$ in $H$. Let $\S$ be the collection of all copies of $C^r_{2\ell}$ generated in this way.
    \begin{claim}
    $$
    |\S|=\Omega\l(e(H)t^{(2\ell-1)(r-1)-1}\r)
    $$
    \end{claim}
    \begin{proof}
        For step (i), there are $e(H)r\,!$ ways to fix an edge with an order. For step (ii), noticing that $t\ge K\ge 4\ell(r-1)$, and that the number of specified vertices is less than $|V(C^r_{2\ell})|=2\ell(r-1)\le t/2$, we know that there are at least $t/2$ ways to specify each of the $2\ell-2$ new vertices. For step (iii), first, there are $(2\ell-1)(r-2)\,!$ ways to fix an order for $f_2,\dots,f_{2\ell}$. Then, by the same argument for step (ii), there are at least $t/2$ ways to specify each of the $(2\ell-1)(r-2)$ new vertices. Therefore, we have $|S|=\Omega(e(H)t^{(2\ell-1)(r-1)-1})$.
    \end{proof}
    \begin{claim}
    $\forall 1\le j\le 2\ell-1$,
    $$
    \Delta_j(\S)=O\l(t^{(2\ell-j)(r-1)-1}\r).
    $$
    \end{claim}
    \begin{proof}
    Let $J$ be a set of $j$ edges of $H$. We want to bound from above the number of ways to specify a copy of $C^r_{2\ell}$ in $\S$ containing $J$. First, we fix which edges of $J$ correspond to which edges of $C^r_{2\ell}$; this may be done in $O(1)$ ways. Note that the number of vertices covered by $j$ edges in $C^r_{2\ell}$ is at least $j(r-1)+1$. Hence there are at most $(2\ell-j)(r-1)-1$ vertices in $C^r_{2\ell}$ that are not specified by the first step. 
    
    Consider a tree $T_{2\ell}$ whose vertices are the edges $e_i, f_i, 1\le i\le 2\ell$, and $e_{i,j}$, $2\le i\le 2\ell$, $1\le j\le r-2$, specified in the algorithm. For $1\le i\le \ell-1$, we let $f_i f_{2\ell+1-i}$ and $f_{2\ell+1-i}f_{i+1}$ be an edge in $T_{2\ell}$; and then for $2\le i\le 2\ell$ and $1\le j\le r-2$, we let $e_{i,j-1}e_{i,j}$ be an edge in $T_{2\ell}$. By definition, if two edges $e$ and $f$ are adjacent in $T$, then they are adjacent in $\Gamma(e\cap f)$. Hence, by the definition of $\Gamma(e\cap f)$, if $e$ is fixed, then the number of ways to specify $f$ is at most $t$. 

    Clearly, if the $C^r_{2\ell}$ has not been completely specified, then there exists an un-specified edge (as a vertex in $T)$ such that one of its neighbor in $T$ is a specified edge. As discussed above, there are at most $t$ ways to specify such an edge. Note that every time we specify such an edge, the number of un-specified vertices will decrease by one. Hence, the $C^r_{2\ell}$ will be specified after at most $(2\ell-j)(r-1)-1$ such edges have been specified. Therefore, we have $\Delta_j(\S)=O(t^{(2\ell-j)(r-1)-1})$.
    \end{proof}
\end{proof}

The balanced supersaturation result for $3$-uniform linear even cycles is slightly different from those when the uniformity is at least $4$. Hence, we state it separately as the following theorem.

\begin{thm}\label{theorem:supersaturation_3}
    For every $\ell\ge 2$, there exist $K, c>0$ such that the following holds for all $n$ sufficiently large and $t\ge K$. Given a $3$-graph $H$ with $n$ vertices and $tn^2$ edges, there exists a collection $\S$ of copies of $C^3_{2\ell}$ in $H$ satisfying 
    $$
    \Delta_j(\S)\le\frac{c|\S|\log n}{tn^2}\l(\max\l\{\l(\frac{t}{\log n}\r)^{-\frac{\ell(4\ell-3)}{4\ell^2-5\ell+2}}n^{-\frac{(\ell-1)(4\ell-3)}{4\ell^2-5\ell+2}},~\l(\frac{t}{\log n}\r)^{-1}n^{-\frac{2\ell-2}{2\ell-1}}\r\}\r)^{j-1},~\forall 1\le j\le 2\ell.
    $$
    In particular, when $t\ge n^{\frac{1}{2(2\ell-1)}}\log n$, 
    $$
    \Delta_j(\S)\le\frac{c|\S|\log n}{tn^2}\l(\l(\frac{t}{\log n}\r)^{-1}n^{-\frac{2\ell-2}{2\ell-1}}\r)^{j-1},~\forall 1\le j\le 2\ell.
    $$
\end{thm}

In the range $t\ge n^{\frac{1}{2(2\ell-1)}}\log n$, a random hypergraph with
$tn^2$ edges shows that Theorem~\ref{theorem:supersaturation_3} is tight up to a polylogarithmic factor. This will be useful in the proof of Theorem~\ref{theorem:supersaturation_ge4}.

\begin{proof}
Let
$$
A=\l(\frac{t}{\log n}\r)^{\frac{(2\ell-1)l}{4\ell^2-5\ell+2}}n^{\frac{(2\ell-1)(\ell-1)}{4\ell^2-5\ell+2}}.
$$  
    We first pick a $3$ partition of $V(H)=V_1\cup V_2\cup V_3$ uniformly at random and let $H'$ be the $3$-partite subgraph of $H$ induced by the partition. Computation on the expected number of edges of $H'$ shows that we can fix an $H'$ with at least $2tn^2/9$ edges. 
     
    we run the following algorithm which partitions $H'$ into subgraphs $F$ and $F_{\tau,a}$, where $0\le a<\log_2A$ and $\tau\in\{\{1,2\},\{1,3\},\{2,3\}\}$.
    \begin{enumerate}
        \item[(i)] Let $H_0=H'$. We construct $H_i$, $i\ge1$, inductively as follows. Given $H_{i-1}$, if there is a $2$-shadow $\sigma$ of $H_{i-1}$ such that $d_{H_{i-1}}(\sigma)<A$, we arbitrarily fix such a $\sigma$ and let $H_i=H_{i-1}\setminus N_{H_{i-1}}(\sigma)$. Then we put all edges of $N_{H_{i-1}}(\sigma)$ into $F_{\tau,a}$ if $\sigma$ contains a vertex of $V_i$ for every $i\in\tau$ and $2^a\le d_{H_{i-1}}(\sigma)< 2^{a+1}$;
        \item[(ii)] If every $2$-shadows $\sigma$ of $H_{i-1}$ has $d_{H_{i-1}}(\sigma)\ge A$, then we let $F=H_{i-1}$. 
    \end{enumerate}
    The proof splits into 2 cases.\\
    
    \noindent\textbf{Case 1: $e(F)\ge e(H')/\log n$.}
    
    In this case, every 2-shadow $\sigma$ of $F$ has $d_F(\sigma)\ge A$. Applying Lemma~\ref{lem: expansion} on $F$ gives us a collection $\S$ of copies of $C^3_{2\ell}$ in $H$ such that
    $$
    |\S|\ge \Omega\l(\frac{tn^2}{\log n} A^{4\ell-3}\r),
    $$
    and, $\forall 1\le j\le 2\ell-1$,
    $$
    \begin{aligned}
        \Delta_j(\S)&\le O\l(A^{4\ell-1-2j}\r)\le O\l(\frac{|\S|\log n}{tn^2A^{4\ell-3}}\cdot A^{4\ell-1-2j}\r)\\
        &\le O\l(\frac{|\S|\log n}{tn^2}\cdot \l(A^{-2}\r)^{j-1}\r)\le O\l(\frac{|\S|\log n}{tn^2}\cdot \l(A^{-\frac{4\ell-3}{2\ell-1}}\r)^{j-1}\r).
    \end{aligned}
    $$

    Note that $\Delta_{2\ell}(\S)=1$. So the bound above also holds for $j=2\ell$.
    Therefore we have, $\forall 1\le j\le 2\ell$,
    $$
    \begin{aligned}
    \Delta_j(\S)&\le O\l(\frac{|\S|\log n}{tn^2}\l(A^{-\frac{4\ell-3}{2\ell-1}}\r)^{j-1}\r)\\
    &= O\l(\frac{|\S|\log n}{tn^2}\l(\l(\frac{t}{\log n}\r)^{-\frac{(4\ell-3)\ell}{4\ell^2-5\ell+2}}n^{-\frac{(4\ell-3)(\ell-1)}{4\ell^2-5\ell+2}}\r)^{j-1}\r)\\
    &\le O\l(\frac{|\S|\log n}{tn^2}\l(\max\l\{\l(\frac{t}{\log n}\r)^{-\frac{\ell(4\ell-3)}{4\ell^2-5\ell+2}}n^{-\frac{(\ell-1)(4\ell-3)}{4\ell^2-5\ell+2}},~\l(\frac{t}{\log n}\r)^{-1}n^{-\frac{2\ell-2}{2\ell-1}}\r\}\r)^{j-1}\r).
    \end{aligned}
    $$
    
    \noindent\textbf{Case 2: $e(F)< e(H')/\log n$.}
    
    Note that the number of edges in $F_{\tau,a}$ is small when $a$ is small. We consider the edges in all $F_{\tau,a}$ with $2^a\le t/9$. By part (i) of the algorithm above, each of these edges contains a 2-shadow with codegree at most $2t/9$. Hence we have,
    $$
    \sum_{\tau\in \binom{[3]}{2}}\sum_{2^a\le t/9} e(F_{\tau,a})\le \binom{n}{2}\frac{2t}{9}\le tn^{2}/9\le e(H')/2.
    $$
    Hence the number of edges in all $F_{\tau,a}$ with $2^a> t/9$ is large,
    $$
    \sum_{\tau\in \binom{[3]}{2}}\sum_{2^a>t/9} e(F_{\tau,a})\ge \Omega(tn^{2}).
    $$

    By the Pigeonhole Principle, there exists a pair of $\tau$ and $a$ with $2^a>t/9$ such that $e(F_{\tau,a})\ge\Omega(tn^2/\log n)$. Arbitrarily fix such a pair of $\tau$ and $a$, let $F'=F_{\tau,a}$ and let $D=2^a$. Since $D>t/9$, $D$ is sufficiently large when $t$ is sufficiently large. Without loss of generality, we can assume $\tau=\{1,2\}$. Let $G$ be the graph on $V_1\cup V_2$ consisting of all $2$-shadows of $F'$ in $V_1\cup V_2$. Note that each 2-shadow has codegree at most $2D$, hence we have $e(G)\ge e(F')/(2D)\ge\Omega (tn^2/(D\log n))$. Let $t'=e(G)/n^{1+1/\ell}$. Then 
    \begin{equation}\label{equation:t'}
        t'\ge\Omega\l(\frac{tn^{1-1/\ell}}{D\log n}\r).
    \end{equation}
    Using $D\le A$, we have
    $$
    t'\ge \Omega\l(\l(\frac{t}{\log n}\r)^{\frac{2(\ell-1)^2}{4\ell^2-5\ell+2}}n^{\frac{2(\ell-1)^3}{\ell(4\ell^2-5\ell+2)}}\r).
    $$
    Hence $t'$ is sufficiently large for us to apply Theorem~\ref{theorem:cycle_supersaturation}, which gives us a collection $\S'$ of copies of $C_{2\ell}$ in $G$ such that 
    \begin{equation}\label{equation:2-shadow}
        \Delta_j(\S')\le O\l(\frac{|\S'|}{t'n^{1+1/\ell}}\l(\max\l\{t'^{-\frac{\ell}{\ell-1}},~t'^{-1}n^{-\frac{\ell-1}{\ell(2\ell-1)}}\r\}\r)^{j-1}\r),~\forall 1\le j\le 2\ell.
    \end{equation}
    
     Note that each $2$-edge of $G$ is contained in at least $D$ $3$-edges of $F'$. Given that $t$ is sufficiently large, we can greedily extend a copy of $C_{2\ell}$ in $G$ into a copy of $C^3_{2\ell}$ in $F'$ in $\Omega(D^{2\ell})$ ways. More formally, let $\S$ be the collection of copies of $C^3_{2\ell}$ in $F'$ with the following property: the 2-shadows of the copy of $C^3_{2\ell}$ in $V_1\cup V_2$ form a $C_{2\ell}$ in $\S'$. Then we have
    \begin{equation}\label{equation:shadow_expansion}
        |\S|\ge\Omega\l(|\S'|D^{2\ell}\r).
    \end{equation}
    
    For any $1\le j\le 2\ell$, given a $j$-tuple of 3-edges of $F'$, by definition, the number of $C_{2\ell}$ in $V_1\cup V_2$ containing the 2-shadows of the $j$ 3-edges is at most $\Delta_j(\S')$. Since each $2$-edge of $G$ is contained in at most $2D$ $3$-edges of $F'$, together with (\ref{equation:t'}), (\ref{equation:2-shadow}) and (\ref{equation:shadow_expansion}) we have, $\forall 1\le j\le 2\ell$,
    $$
    \begin{aligned}
    \Delta_j(\S)&\le O\l(\Delta_j(\S')D^{2\ell-j}\r)\\
    &\le O\l(\frac{|\S|}{t'n^{1+1/\ell}}\l(\max\l\{t'^{-\frac{\ell}{\ell-1}},~t'^{-1}n^{-\frac{\ell-1}{\ell(2\ell-1)}}\r\}\r)^{j-1}D^{-j}\r)\\
    &\le O\l(\frac{|\S|\log n}{tn^{2}}\l(\max\l\{\l(\frac{t}{\log n}\r)^{-\frac{\ell}{\ell-1}}n^{-1}D^{\frac{1}{\ell-1}},~\l(\frac{t}{\log n}\r)^{-1}n^{-\frac{2\ell-2}{2\ell-1}}\r\}\r)^{j-1}\r).
    \end{aligned}
    $$

    Finally, note that
    $$
D\le A=\l(\frac{t}{\log n}\r)^{\frac{(2\ell-1)l}{4\ell^2-5\ell+2}}n^{\frac{(2\ell-1)(\ell-1)}{4\ell^2-5\ell+2}},
$$  
hence we have
    $$
         \Delta_j(\S)\le O\l(\frac{|\S|\log n}{tn^2}\l(\max\l\{\l(\frac{t}{\log n}\r)^{-\frac{\ell(4\ell-3)}{4\ell^2-5\ell+2}}n^{-\frac{(\ell-1)(4\ell-3)}{4\ell^2-5\ell+2}},~\l(\frac{t}{\log n}\r)^{-1}n^{-\frac{2\ell-2}{2\ell-1}}\r\}\r)^{j-1}\r).
    $$
\end{proof}

Next, we prove the balanced supersaturation theorems for $C^r_{2\ell}$ when $r\ge 4$.

\begin{thm}\label{theorem:supersaturation_ge4}
For every $r\ge 4$ and $\ell\ge 2$, there exist $K, c>0$ such that the following holds for all $n$ sufficiently large and $t\ge K$. Given an $r$-graph $H$ with $n$ vertices and $tn^{r-1}$ edges, there exists a collection $\S$ of copies of $C^r_{2\ell}$ in $H$ satisfying
$$    
\Delta_j(\S)\le\frac{c|\S|(\log n)^{r-2}}{tn^{r-1}}\l(\l(\frac{t}{(\log n)^{r-2}}\r)^{-1}n^{-r+2+\frac{1}{2\ell-1}}\r)^{j-1},~\forall 1\le j\le 2\ell.     
$$
\end{thm}

A random hypergraph with $tn^{r-1}$ edges shows that Theorem~\ref{theorem:supersaturation_ge4} is tight up to a polylogarithmic factor. Recall that Theorem~\ref{theorem:supersaturation_3} is tight when $t$ is large enough (that is, the given 3-graph is dense enough). The reason why one can prove a tight balanced supersaturation for all $t\ge \Omega(1)$ when $r=4$ is that when the ``typical'' codegree is ``small'', the 3-shadows will be dense enough to miraculously land in the tight range of Theorem~\ref{theorem:supersaturation_3} (see the comment after the statement of Theorem~\ref{theorem:supersaturation_3}). After proving a tight result for $r=4$, one can then obtain a tight result for all $r\ge 4$ by induction.

\begin{proof}
    We use induction on $r$. When $r=4$, we first pick a $4$-partition of $V(H)=V_1\cup V_2\cup V_3\cup V_4$ uniformly at random and let $H'$ be the $4$-partite subgraph of $H$ induced by the partition. Computation on the expected number of edges of $H'$ shows that we can fix an $H'$ with at least $3tn^3/32$ edges. Let 
    $$A_4=\l(\frac{t}{(\log n)^2}\r)^{\frac{2\ell-1}{6\ell-4}}n^{\frac{4\ell-3}{6\ell-4}}.$$
    We now describe an algorithm which partitions $H'$ into subgraphs $F$ and $F_{\tau,a}$, where $0\le a<\log_2A_4$ and $\tau\in\{\{1,2,3\},\{1,2,4\},\{1,3,4\},\{2,3,4\}\}$.
    \begin{enumerate}
        \item[(i)] Let $H_0=H'$. We construct $H_i$, $i\ge1$, inductively as follows. Given $H_{i-1}$, if there is a $3$-shadow $\sigma$ of $H_{i-1}$ such that $d_{H_{i-1}}(\sigma)<A_4$, we fix such a $\sigma$ and let $H_i=H_{i-1}\setminus N_{H_{i-1}}(\sigma)$. Then we put all edges of $N_{H_{i-1}}(\sigma)$ into $F_{\tau,a}$ if $\sigma$ contains a vertex of $V_i$ for every $i\in\tau$ and $2^a\le d_{H_{i-1}}(\sigma)< 2^{a+1}$;
        \item[(ii)] If every $3$-shadow $\sigma$ of $H_{i-1}$ has $d_{H_{i-1}}(\sigma)\ge A_4$, then we let $F=H_{i-1}$. 
    \end{enumerate}
    The proof splits into 2 cases.\\

    \noindent\textbf{Case 1: $e(F)\ge e(H')/(\log n)^2$.}
    
    In this case, every 3-shadow $\sigma$ of $F$ has $d_F(\sigma)\ge A_4$. Applying Lemma~\ref{lem: expansion} on $F$ gives us a collection $\S$ of copies of $C^4_{2\ell}$ in $F$ satisfying
    $$
    |\S|\ge\Omega\l(\frac{tn^3}{(\log n)^2}A_4^{6\ell-4}\r),~\Delta_j(\S)\le O\l(A_4^{6l-3j-1}\r),~\forall 1\le j\le 2\ell-1.
    $$
    Also note that $\Delta_{2\ell}(\S)=1$. Hence we have, $\forall 1\le j\le 2\ell$,
    $$
    \begin{aligned}
        \Delta_j(\S)&\le O\l(\frac{|S|(\log n)^2}{tn^3}\l(A_4^{-\frac{6\ell-4}{2\ell-1}}\r)^{j-1}\r)\\
        &=O\l(\frac{|S|(\log n)^2}{tn^3}\l(\l(\frac{t}{(\log n)^2}\r)^{-1}n^{-\frac{4\ell-3}{2\ell-1}}\r)^{j-1}\r).
    \end{aligned}
    $$
    \medskip
    
    \noindent\textbf{Case 2: $e(F)<e(H')/(\log n)^2$.}

    Note that the number of edges in $F_{\tau,a}$ is small when $a$ is small. We consider the number of edges in all $F_{\tau,a}$ with $2^a\le 3t/16$,
    $$
    \sum_{\tau\in \binom{[4]}{3}}\sum_{2^a\le 3t/16} e(F_{\tau,a})\le \binom{n}{3}\frac{3t}{8}\le \frac{tn^{3}}{16}\le 2e(H')/3.
    $$
    Hence the number of edges in all $F_{\tau,a}$ with $2^a> 3t/16$ is large,
    $$
    \sum_{\tau\in \binom{[4]}{3}}\sum_{2^a>3t/16} e(F_{\tau,a})\ge \Omega(tn^{3}).
    $$

    By the Pigeonhole Principle, there exists a pair of $\tau$ and $a$ with $2^a> 3t/16$ such that $e(F_{\tau,a})\ge\Omega\l(tn^3/\log n\r)$. Arbitrarily fix such a pair of $\tau$ and $a$, let $F'=F_{\tau,a}$ and let $D=2^a$. Since $D>3t/16$, $D$ is sufficiently large when $t$ is sufficiently large. Without loss of generality, we can assume $\tau=\{1,2,3\}$. Let $G$ be the 3-graph on $V_1\cup V_2\cup V_3$ consisting of all $3$-shadows of $F'$ in $V_1\cup V_2\cup V_3$. Then $e(G)\ge e(F')/(2D)\ge\Omega (tn^3/(D\log n))$. Let $t'=e(G)/n^2\ge\Omega(tn/(D\log n))$. Note that when $n$ is sufficiently large, 
    $$
    t'\ge\Omega\l(\frac{tn}{A_4\log n}\r)=\Omega\l(\log n\l(\frac{t}{(\log n)^2}\r)^\frac{4\ell-3}{6\ell-4}n^{\frac{2\ell-1}{6\ell-4}}\r)\ge  n^{\frac{1}{2(2\ell-1)}}\log n.
    $$ 
    Hence by Theorem~\ref{theorem:supersaturation_3}, there exists a collection $\S'$ of copies of $C^3_{2\ell}$ in $G$ such that 
    $$
    \Delta_j(\S')\le O\l(\frac{|\S'|\log n}{t'n^2}\l(\l(\frac{t'}{\log n}\r)^{-1}n^{-\frac{2\ell-2}{2\ell-1}}\r)^{j-1}\r),~\forall 1\le j\le 2\ell.
    $$
    Note that each $3$-edge of $G$ is contained in at least $D$ $4$-edges of $F'$. Given that $t$ is sufficiently large, we can greedily extend a copy of $C^3_{2\ell}$ in $G$ into a copy of $C^4_{2\ell}$ in $F'$ in $\Omega(D^{2\ell})$ ways. More formally, let $\S$ be the collection of copies of $C^4_{2\ell}$ in $F'$ with the following property: the $3$-shadows of the copy of $C^4_{2\ell}$ in $V_1\cup V_2\cup V_3$ form a $C^3_{2\ell}$ in $\S'$. Then
    $$
    |\S|\ge\Omega\l(|\S'|D^{2\ell}\r).
    $$
    Since each $3$-shadow of $F'$ in $V_1\cup V_2\cup V_3$ is contained in at most $2D$ edges of $F'$, we have, $\forall 1\le j\le 2\ell$,
    $$
    \begin{aligned}
        \Delta_j(\S)&\le O\l(\Delta_j(\S')D^{2\ell-j}\r)\\
        &\le O\l(\frac{|\S|(\log n)^2}{tn^3}D^{-j+1}\l(\l(\frac{t}{(\log n)^2D}\r)^{-1}n^{-\frac{4\ell-3}{2\ell-1}}\r)^{j-1}\r)\\
        &=O\l(\frac{|\S|(\log n)^2}{tn^3}\l(\l(\frac{t}{(\log n)^2}\r)^{-1}n^{-\frac{4\ell-3}{2\ell-1}}\r)^{j-1}\r).
    \end{aligned}
    $$   
    This completes the proof for $r=4$.

    \medskip

    When $r\ge5$, we simply repeat the argument above with
    $$
    A_r=\l(\frac{t}{(\log n)^{r-2}}\r)^{\frac{2\ell-1}{(2\ell-1)r-2\ell}}n^{\frac{(2\ell-1)r-4\ell+1}{(2\ell-1)r-2\ell}}.
    $$

We pick an $r$-partition of $V(H)=V_1\cup\dots\cup V_r$ uniformly at random and let $H'$ be the $r$-partite subgraph of $H$ induced by the partition. Computation on the expected number of edges of $H'$ shows that we can fix an $H'$ with at least $r\,!\,tn^{r-1}/r^r$ edges. 

We now describe an algorithm which partitions $H'$ into subgraphs $F$ and $F_{\tau,a}$, where $0\le a<\log_2A_r$ and $\tau\in\binom{[r]}{r-1}$.
    \begin{enumerate}
        \item[(i)] Let $H_0=H'$. We construct $H_i$, $i\ge1$, inductively as follows. Given $H_{i-1}$, if there is an $(r-1)$-shadow $\sigma$ of $H_{i-1}$ such that $d_{H_{i-1}}(\sigma)<A_r$, we fix such a $\sigma$ and let $H_i=H_{i-1}\setminus N_{H_{i-1}}(\sigma)$. Then we put all edges of $N_{H_{i-1}}(\sigma)$ into $F_{\tau,a}$ if $\sigma$ contains a vertex of $V_i$ for every $i\in\tau$ and $2^a\le d_{H_{i-1}}(\sigma)< 2^{a+1}$;
        \item[(ii)] If every $(r-1)$-shadow $\sigma$ of $H_{i-1}$ has $d_{H_{i-1}}(\sigma)\ge A_r$, then we let $F=H_{i-1}$. 
    \end{enumerate}
    The proof splits into 2 cases.\\

    \noindent\textbf{Case 1: $e(F)\ge e(H')/(\log n)^{r-2}$.}
    
    In this case, every $(r-1)$-shadow $\sigma$ of $F$ has $d_F(\sigma)\ge A_r$. Applying Lemma~\ref{lem: expansion} on $F$ gives us a collection $\S$ of copies of $C^r_{2\ell}$ in $F$ satisfying
    $$
    |\S|\ge\Omega\l(\frac{tn^{r-1}}{(\log n)^{r-2}}A_r^{(2\ell-1)r-2\ell}\r),~\Delta_j(\S)\le O\l(A_r^{(2\ell-1)r-2\ell-(j-1)(r-1)}\r),~\forall 1\le j\le 2\ell-1.
    $$
    Also note that $\Delta_{2\ell}(\S)=1$. Hence we have, $\forall 1\le j\le 2\ell$,
    $$
    \begin{aligned}
        \Delta_j(\S)&\le O\l(\frac{|S|(\log n)^{r-2}}{tn^{r-1}}\l(A_r^{-\frac{(2\ell-1)r-2\ell}{2\ell-1}}\r)^{j-1}\r)\\
        &=O\l(\frac{|S|(\log n)^{r-2}}{tn^{r-1}}\l(\l(\frac{t}{(\log n)^{r-2}}\r)^{-1}n^{-r+2+\frac{1}{2\ell-1}}\r)^{j-1}\r).
    \end{aligned}
    $$
    \medskip
    
    \noindent\textbf{Case 2: $e(F)<e(H')/(\log n)^{r-2}$.}

    Note that the number of edges in $F_{\tau,a}$ is small when $a$ is small. We consider the number of edges in all $F_{\tau,a}$ with $2^a\le \frac{r\,!}{4r^r}t$,
    $$
    \sum_{\tau\in \binom{[r]}{r-1}}\sum_{2^a\le \frac{r\,!}{4r^r}t} e(F_{\tau,a})\le \binom{n}{r-1}\frac{r\,!}{2r^r}t\le \frac{r}{2r^r} tn^{r-1}\le e(H')/2.
    $$
    Hence the number of edges in all $F_{\tau,a}$ with $2^a> \frac{r\,!}{4r^r}t$ is large, 
    $$
    \sum_{\tau\in \binom{[r]}{r-1}}\sum_{2^a>\frac{r\,!}{4r^r}t} e(F_{\tau,a})\ge \Omega(tn^{r-1}).
    $$

    By the Pigeonhole Principle, there exists a pair of $\tau$ and $a$ with $2^a> \frac{r\,!}{4r^r}t$ such that $e(F_{\tau,a})\ge\Omega\l(tn^{r-1}/\log n\r)$. Arbitrarily fix such a pair of $\tau$ and $a$, let $F'=F_{\tau,a}$ and let $D=2^a$. Since $D>\frac{r\,!}{4r^r}t$, $D$ is sufficiently large when $t$ is sufficiently large. Without loss of generality, we can assume $\tau=\{1,\dots,r-1\}$. Let $G$ be the $(r-1)$-graph on $V_1\cup \dots\cup V_{r-1}$ consisting of all $(r-1)$-shadows of $F'$ in $V_1\cup \dots\cup V_{r-1}$. Then $e(G)\ge e(F')/(2D)\ge\Omega (tn^{r-1}/(D\log n))$. Let $t'=e(G)/n^{r-2}$. Then $t'\ge\Omega(tn/(D\log n))$. Since $2D\le A_r\le o(tn/\log n)$, $t'$ is sufficiently large given that $n$ is sufficiently large.
    
    Hence by inductive hypothesis, there exists a collection $\S'$ of copies of $C^{r-1}_{2\ell}$ in $G$ such that 
    $$
    \Delta_j(\S')\le O\l(\frac{|\S'|(\log n)^{r-3}}{t'n^{r-2}}\l(\l(\frac{t'}{(\log n)^{r-3}}\r)^{-1}n^{-r+3+\frac{1}{2\ell-1}}\r)^{j-1}\r),~\forall 1\le j\le 2\ell.
    $$
    Note that each $(r-1)$-edge of $G$ is contained in at least $D$ $r$-edges of $F'$. Given that $t$ is sufficiently large, we can greedily extend a copy of $C^{r-1}_{2\ell}$ in $G$ into a copy of $C^r_{2\ell}$ in $F'$ in $\Omega(D^{2\ell})$ ways. More formally, let $\S$ be the collection of copies of $C^r_{2\ell}$ in $F'$ with the following property: the $(r-1)$-shadows of the copy of $C^r_{2\ell}$ in $V_1\cup \dots\cup V_{r-1}$ form a $C^{r-1}_{2\ell}$ in $\S'$. Then
    $$
    |\S|\ge\Omega\l(|\S'|D^{2\ell}\r).
    $$
    Since each $(r-1)$-shadow of $F'$ in $V_1\cup\dots\cup V_{r-1}$ is contained in at most $2D$ edges of $F'$, we have, $\forall 1\le j\le 2\ell$,
    $$
    \begin{aligned}
        \Delta_j(\S)&\le O\l(\Delta_j(\S')D^{2\ell-j}\r)\\
        &\le O\l(\frac{|\S|(\log n)^{r-2}}{tn^{r-1}}D^{-j+1}\l(\l(\frac{tn}{(\log n)^{r-2}D}\r)^{-1}n^{-r+3+\frac{1}{2\ell-1}}\r)^{j-1}\r)\\
        &=O\l(\frac{|\S|(\log n)^{r-2}}{tn^{r-1}}\l(\l(\frac{t}{(\log )^{r-2}}\r)^{-1}n^{-r+2+\frac{1}{2\ell-1}}\r)^{j-1}\r).
    \end{aligned}
    $$   

%%%%%%%%%%%%%%%%%%%%%%%%%%%%%%%%%%%%%%%%%%%%%%%%%%%%%%%%%%%%%%%
\end{proof}

%%%%%%%%%%%%%%%%%%%%%%%%%%%%%%%%%%%%%%%%%%%%%%%%%%%%%%%%%%%%%%%%%%%%%%%%%%%%%%%%%%%%%%%%%%%%%%%%%%%%%%%%%%%%%%%%%%%%%%%%%%%%%%%%%%%%%%%%%%%%%%%%%%%%%%%%%%%%%%%%%%%%%%%%%%%%%%%%%%%%%%%%%%%%%%%%%%%%%%%%%%%%%%%%%%%%%%%%%%%%%%%%%%%%%%%%%%%%%%%%%%%%%%%%%%%%%%%%%%%%%%%%%%%%%%%%%%%%%%%%%%%%%%%%%%%%%%%%%%%%%%%%%%%%%%%%%%%%%%%%%%%%%%%%%%%%%%%%%%%%%%%%%%%%%%%%%%%%%%%%%%%%%%%%%%%%%%%%%%%%%%%%%%%%%%%%%%%%%%%%%%%%%%%%%%%%%%%%%%%%%%%%%%%%%%%%%%%%%%%%%%%%%%%%%%%%%%%%%%%%%%%%%%%%%%%%%%%%%%%%%%%%%%%%%%%%%%%

\section{Container theorems}\label{section:container}

We make use of the following simplified version of the Containers Theorem by Balogh, Morris and Samotij~\cite{balogh2018method}:
\begin{thm}[\cite{balogh2018method}]\label{theorem:container}
For every $r\ge 2$, there exists a constant $\epsilon>0$ such that the following holds.
Let $\S$ be an $r$-graph such that
    \begin{equation}\label{equation:container}
        \Delta_j(\S)\leq \l(\frac{B}{v(\S)}\r)^{j-1}\frac{e(\S)}{L},~~\forall 1\le j\le r, 
    \end{equation}
    for some integers $B, L>0$.
     Then there exists a collection $\mathcal{C}$ of at most
     $$
     \exp\l(\frac{\log\l(\frac{v(\S)}{B}\r)B}{\epsilon}\r)
     $$
     subsets of $V(\S)$ such that:
    \begin{enumerate}
        \item[(a)] for every independent set $I$ of $\S$, there exists $C\in \mathcal{C}$ such that $I\subset C$;
        \item[(b)] for every $C\in\mathcal{C}$, $|C|\leq v(\S)-\epsilon L$.
    \end{enumerate}
\end{thm}

Consider a hypergraph $\S$ on $E(H)$ whose edges are copies of $C^r_{2\ell}$ guaranteed by Theorem~\ref{theorem:supersaturation_3} or Theorem~\ref{theorem:supersaturation_ge4}. By applying Theorem~\ref{theorem:container} on $\S$ we have the following theorems.

\begin{thm}\label{theorem:container1_3}
    For every $\ell\ge2$, there exist $K,\epsilon>0$ such that the following holds for all sufficiently large $n$ and every $t\ge K$. Given a $3$-graph $H$ with $n$ vertices and $tn^{2}$ edges there exists a collection $\C$ of at most 
    $$
    \exp\l(\frac{(\log n)^2}{\epsilon}\max\l\{n^{\frac{4\ell^2-3\ell+1}{4\ell^2-5\ell+2}}\l(\frac{\log n}{t}\r)^{\frac{2\ell-2}{4\ell^2-5\ell+2}},n^{\frac{2\ell}{2\ell-1}}\r\}\r)
    $$
    subgraphs of $H$ such that
    \begin{itemize}
        \item[(a)] every $C^3_{2\ell}$-free subgraph of $H$ is a subgraph of some $G\in\C$;
        \item[(b)] $e(G)\le \l(1-\frac{\epsilon}{\log n}\r)tn^2$, for each $G\in\C$. 
    \end{itemize}
    % In particular, when $t\ge n^{\frac{1}{4\ell-2}}\log n$, or in other words, when $|H|\ge n^{\frac{8l-3}{4\ell-2}}\log n$,
    % $$|\C|\le \exp\l(\frac{(\log n)^2n^{\frac{2\ell}{2\ell-1}}}{\epsilon}\r)$$
\end{thm}

\begin{proof}
    Given that $K$ is sufficiently large, we can apply Theorem~\ref{theorem:supersaturation_3} on $H$. Let $c$ be the constant and let $\S$ be the collection of copies of $C^3_{2\ell}$ in $H$ guaranteed by Theorem~\ref{theorem:supersaturation_3}. We view $\S$ as a hypergraph on $E(H)$ and apply on it Theorem~\ref{theorem:container} with parameters
    $$
    L=\frac{tn^2}{c\log n},~\text{and}~B=e(H)\max\l\{\l(\frac{t}{\log n}\r)^{-\frac{\ell(4\ell-3)}{4\ell^2-5\ell+2}}n^{-\frac{(\ell-1)(4\ell-3)}{4\ell^2-5\ell+2}},~\l(\frac{t}{\log n}\r)^{-1}n^{-\frac{2\ell-2}{2\ell-1}}\r\}.
    $$
    This gives a collection $\C$ of subgraphs of $H$ and $\epsilon=\epsilon(\ell)$ such that
    \begin{equation}\label{equation:num_container_3}
        |\C|\le \exp\l(\frac{(\log n)^2}{\epsilon}\max\l\{\l(\frac{t}{\log n}\r)^{-\frac{2\ell-2}{4\ell^2-5\ell+2}}n^{\frac{4\ell^2-3\ell+1}{4\ell^2-5\ell+2}},~n^{\frac{2\ell}{2\ell-1}}\r\}\r);
    \end{equation}
    every $C^3_{2\ell}$-free subgraph of $H$ is a subgraph of some $G\in\C$; and
    \begin{equation}\label{equation:size_container_3}
        e(G)\le \l(1-\frac{\epsilon}{\log n}\r)tn^2
    \end{equation}
    for each $G\in\C$. Note that, for simplicity, constants factors in (\ref{equation:num_container_3}) and (\ref{equation:size_container_3}) are absorbed by $\epsilon$ which is sufficiently small.
\end{proof}

\begin{thm}\label{theorem:container1_ge4}
    For every $r\ge 4$ and $\ell\ge2$, there exist $K,\epsilon>0$ such that the following holds for all sufficiently large $n$ and every $t\ge K$. Given an $r$-graph $H$ with $n$ vertices and $tn^{r-1}$ edges there exists a collection $\C$ of at most 
    $$
    \exp\l(\frac{(\log n)^{r-1}n^{\frac{2\ell}{2\ell-1}}}{\epsilon}\r)
    $$
    subgraphs of $H$ such that
    \begin{itemize}
        \item[(a)] every $C^r_{2\ell}$-free subgraph of $H$ is a subgraph of some $G\in\C$;
        \item[(b)] $e(G)\le \l(1-\frac{\epsilon}{(\log n)^{r-2}}\r)tn^{r-1}$, for each $G\in\C$. 
    \end{itemize}
\end{thm}

\begin{proof}
   Given that $K$ is sufficiently large, we can apply Theorem~\ref{theorem:supersaturation_ge4} on $H$. Let $c$ be the constant and let $\S$ be the collection of copies of $C^r_{2\ell}$ in $H$ guaranteed by Theorem~\ref{theorem:supersaturation_ge4}. We view $\S$ as a hypergraph on $E(H)$ and apply on it Theorem~\ref{theorem:container} with parameters
    $$
    L=\frac{tn^{r-1}}{c(\log n)^{r-2}},~\text{and}~B=e(H)\l(\frac{t}{(\log n)^{r-2}}\r)^{-1}n^{-r+2+\frac{1}{2\ell-1}}.
    $$
    This gives a collection $\C$ of subgraphs of $H$ and $\epsilon=\epsilon(\ell)$ such that
    \begin{equation}\label{equation:num_container_ge4}
            |\C|\le \exp\l(\frac{(\log n)^{r-1}n^{\frac{2\ell}{2\ell-1}}}{\epsilon}\r);
    \end{equation}
    every $C^r_{2\ell}$-free subgraph of $H$ is a subgraph of some $G\in\C$; and
    \begin{equation}\label{equation:size_container_ge4}
        e(G)\le \l(1-\frac{\epsilon}{(\log n)^{r-2}}\r)tn^{r-1}    
    \end{equation}
    for each $G\in\C$. Note that, for simplicity, constants factors in (\ref{equation:num_container_ge4}) and (\ref{equation:size_container_ge4}) are absorbed by $\epsilon$ which is sufficiently small.
\end{proof}

We obtain the following theorem by applying Theorem~\ref{theorem:container1_3} repeatedly.

\begin{thm}\label{theorem:container2_3}
    For every $\ell\ge2$, there exist $c, K>0$ such that for all sufficiently large integers $n$ and $t\ge K$, there exists a collection $\C$ of $3$-graphs on $[n]$ such that 
    \begin{enumerate}
        \item[(a)]
        $$
        |\C|\le \exp\l(c\cdot\max\l\{(\log n)^{3+\frac{2\ell-2}{4\ell^2-5\ell+2}}n^{\frac{4\ell^2-3\ell+1}{4\ell^2-5\ell+2}}t^{-\frac{2\ell-2}{4\ell^2-5\ell+2}},~(\log n)^4 n^{\frac{2\ell}{2\ell-1}}\r\}\r);
        $$
        \item[(b)] $e(G)\le tn^{2}$, $\forall G\in\C$;
        \item[(c)] every $C^3_{2\ell}$-free $3$-graph on $[n]$ is a subgraph of some $G\in\C$.
    \end{enumerate}
    In particular, when $t\le n^{\frac{1}{4\ell-2}}\l(\log n\r)^{-\frac{4\ell^2-7\ell+4}{2\ell-2}}$, 
    $$
    |\C|\le \exp\l(c\cdot(\log n)^{3+\frac{2\ell-2}{4\ell^2-5\ell+2}}n^{\frac{4\ell^2-3\ell+1}{4\ell^2-5\ell+2}}t^{-\frac{2\ell-2}{4\ell^2-5\ell+2}}\r).
    $$
\end{thm}

\begin{proof}
    Take $K=K(\ell)$ and $\epsilon=\epsilon(\ell)$ from Theorem~\ref{theorem:container1_3}. Let $H_0=K^3_n$ and let $\C_0=\{H_0\}$. Define $t_0=\binom{n}{3}/n^2$, $t_i=\exp\l(-\epsilon/\log n\r)t_{i-1}$ for all $i\ge 1$. Let $m$ be the smallest integer such that $t_m\le t$. By applying Theorem~\ref{theorem:container1_3} on $H_0=K^3_n$, we have a collection $\C_1$ of subgraphs of $K^3_n$ such that every $C^3_{2\ell}$-free subgraph of $K^3_n$ is a subgraph of some $G\in\C_1$, and $e(G)\le \exp(-\epsilon/\log n)t_0n^2=t_1n^2$ for each $G\in\C_1$.

    For $1\le i\le m-1$, given $\C_i$ such that $e(G)\le t_in^2$ for each $G\in \C_i$, we construct $\C_{i+1}$ as follows: for every $G\in\C_i$, if $e(G)\le t_{i+1}n^2$, then put $G$ into $\C_{i+1}$; if $e(G)>t_{i+1}n^2$, we apply Theorem~\ref{theorem:container1_3} on $G$ to obtain a collection $\C'$ of subgraphs of $G$ such that $e(G')\le t_{i+1}n^2$ for each $G'\in\C'$, and then put every element in $\C'$ into $\C_{i+1}$. Let $\C=\C_m$. Clearly, every $C^3_{2\ell}$-free subgraph of $K^3_n$ is a subgraph of some $G\in\C$, and $e(G)\le tn^2$ for each $G\in\C$. Also, for $0\le i\le m-1$, we have
    $$
    \begin{aligned}
    \frac{|\C_{i+1}|}{|\C_i|}&\le \exp\l(\frac{(\log n)^2}{\epsilon}\max\l\{n^{\frac{4\ell^2-3\ell+1}{4\ell^2-5\ell+2}}\l(\frac{\log n}{t_{i+1}}\r)^{\frac{2\ell-2}{4\ell^2-5\ell+2}},n^{\frac{2\ell}{2\ell-1}}\r\}\r)\\
    &\le \exp\l(\frac{(\log n)^2}{\epsilon}\l(n^{\frac{4\ell^2-3\ell+1}{4\ell^2-5\ell+2}}\l(\frac{\log n}{t_{i+1}}\r)^{\frac{2\ell-2}{4\ell^2-5\ell+2}}+n^{\frac{2\ell}{2\ell-1}}\r)\r).
    \end{aligned}
    $$
    It is not hard to check that $m\le (\log n)^2/\epsilon$. Observe that $\l\{t_{i}^{-\frac{2\ell-2}{4\ell^2-5\ell+2}}\r\}_{0\le i\le m-1}$ is an increasing geometric sequence with common ratio $\exp\l(\frac{(2\ell-2)\epsilon}{(4\ell^2-5\ell+2)\log n}\r)$, so its sum
    $$
    \sum_{i=0}^{m-1}t_i^{-\frac{2\ell-2}{4\ell^2-5\ell+2}}\le O\l(t_{m-1}^{-\frac{2\ell-2}{4\ell^2-5\ell+2}}\cdot\frac{1}{1-\exp\l(-\frac{(2\ell-2)\epsilon}{(4\ell^2-5\ell+2)\log n}\r)}\r)\le O\l(t^{-\frac{2\ell-2}{4\ell^2-5\ell+2}}\log n\r).
    $$
    
    Hence there exists a constant $c>0$ such that
    $$
    \begin{aligned}
        |\C|&\le\prod_{i=0}^{m-1}\frac{|\C_{i+1}|}{|\C_i|}\\
         &\le\exp\l(\frac{(\log n)^2}{\epsilon}\l(n^{\frac{4\ell^2-3\ell+1}{4\ell^2-5\ell+2}}\l({\log n}\r)^{\frac{2\ell-2}{4\ell^2-5\ell+2}}\sum_{i=0}^{m-1}t_i^{-\frac{2\ell-2}{4\ell^2-5\ell+2}}+mn^{\frac{2\ell}{2\ell-1}}\r)\r)\\
         &\le \exp\l(c\cdot\max\l\{(\log n)^{3+\frac{2\ell-2}{4\ell^2-5\ell+2}}n^{\frac{4\ell^2-3\ell+1}{4\ell^2-5\ell+2}}t^{-\frac{2\ell-2}{4\ell^2-5\ell+2}},~(\log n)^4 n^{\frac{2\ell}{2\ell-1}}\r\}\r).
    \end{aligned}
    $$
\end{proof}

The following theorem is obtained by applying Theorem~\ref{theorem:container1_ge4} repeatedly.

\begin{thm}\label{theorem:container2_ge4}
    For every $\ell\ge2$ and $r\ge4$, there exist $c, K>0$ such that for all sufficiently large $n$ and every $t\ge K$, there exists a collection $\C$ of $r$-graphs on $[n]$ such that 
    \begin{enumerate}
        \item[(a)]
        $$
        |\C|\le \exp\l(c(\log n)^{2r-2}n^{\frac{2\ell}{2\ell-1}}\r)
        $$
        \item[(b)] $e(G)\le tn^{r-1}$, $\forall G\in\C$;
        \item[(c)] every $C^r_{2\ell}$-free $r$-graph on $[n]$ is a subgraph of some $G\in\C$.
    \end{enumerate}
\end{thm}

\begin{proof}
    Take $K=K(\ell)$ and $\epsilon=\epsilon(\ell)$ from Theorem~\ref{theorem:container1_ge4}. Let $H_0=K^r_n$, $\C_0=\{H_0\}$. Define $t_0=\binom{n}{r}/n^{r-1}$, $t_i=\exp(-\epsilon/(\log n)^{r-2})t_{i-1}$ for all $i\ge 1$. Let $m$ be the smallest integer such that $t_m\le t$. By applying Theorem~\ref{theorem:container1_ge4} on $H_0=K^r_n$, we have a collection $\C_1$ of subgraphs of $K^r_n$ such that every $C^r_{2\ell}$-free subgraph of $K^r_n$ is a subgraph of some $G\in\C_1$, and $e(G)\le \exp(-\epsilon/(\log n)^{r-2})t_0n^{r-1}=t_1n^{r-1}$ for each $G\in\C_1$.

    For $1\le i\le m-1$, given $\C_i$ such that $e(G)\le t_in^{r-1}$ for each $G\in \C_i$, we construct $\C_{i+1}$ as follows: for every $G\in\C_i$, if $e(G)\le t_{i+1}n^{r-1}$, then we put $G$ into $\C_{i+1}$; if $e(G)>t_{i+1}n^{r-1}$, we apply Theorem~\ref{theorem:container1_ge4} on $G$ to obtain a collection $\C'$ of subgraphs of $G$ such that $e(G')\le t_{i+1}n^{r-1}$ for each $G'\in\C'$, and then put every element in $\C'$ into $\C_{i+1}$. 
    
    Let $\C=\C_m$. Clearly, every $C^r_{2\ell}$-free subgraph of $K^r_n$ is a subgraph of some $G\in\C$, and $e(G)\le tn^{r-1}$ for each $G\in\C$. Also, for $0\le i\le m-1$, we have
    $$
    \frac{|\C_{i+1}|}{|\C_i|}\le\exp\l(\frac{(\log n)^{r-1}n^{\frac{2\ell}{2\ell-1}}}{\epsilon}\r).
    $$
     It is not hard to check that $m\le (\log n)^{r-1}/\epsilon$. Hence, there exists a constant $c>0$ such that
    $$
    |\C|=\prod_{i=1}^{m-1}\frac{|\C_{i+1}|}{|\C_i|}\le \exp\l(c{(\log n)^{2r-2}n^{\frac{2\ell}{2\ell-1}}}\r).
    $$
\end{proof}

%%%%%%%%%%%%%%%%%%%%%%%%%%%%%%%%%%%%%%%%%%%%%%%%%%%%%%%%%%%%%%%%%%%%%%%%%%%%%%%%%%%%%%%%%%%%%%%%%%%%%%%%%%%%%%%%%%%%%%%%%%%%%%%%%%%%%%%%%%%%%%%%%%%%%%%%%%%%%%%%%%%%%%%%%%%%%%%%%%%%%%%%%%%%%%%%%%%%%%%%%%%%%%%%%%%%%%%%%%%%%%%%%%%%%%%%%%%%%%%%%%%%%%%%%%%%%%%%%%%%%%%%%%%%%%%%%%%%%%%%%%%%%%%%%%%%%%%%%%%%%%%%%%%%%%%%%%%%%%%%%%%%%%%%%%%%%%%%%%%%%%%%%%%%%%%%%%%%%%%%%%%%%%%%%%%%%%%%%%%%%%%%%%%%%%%%%%%%%%%%%%%%%%%%%%%%%%%%%%%%%%%%%%%%%%%%%%%%%%%%%%%%%%%%%%%%%%%%%%%%%%%%%%%%%%%%%%%%%%%%%%%%%%%%%%%%%%%%

\section{Random Tur\'an Theorems for linear even cycles}\label{section:upper}
In this section, we prove Theorem~\ref{theorem:main_ge4} and Theorem~\ref{theorem:main_3}.

% \setcounter{thmB}{\getrefnumber{theorem:main_ge4}}
% \addtocounter{thmB}{-1}

% \begin{thmB}
%     For every $r\ge 4$ and $\ell\ge2$, there exists $c>0$ such that if $p\ge (\log n)^{2r-2}n^{-r+2+\frac{1}{2\ell-1}}$, then a.a.s.
%     $$
%     \ex(G^r_{n,p}, C^r_{2\ell})\le cpn^{r-1}.
%     $$
% \end{thmB}

\begin{proof}[Proof of Theorem~\ref{theorem:main_ge4}]
Recall that we wish to prove that for every $r\ge 4$ and $\ell\ge2$, there exists $c>0$ such that if $p\ge (\log n)^{2r-2}n^{-r+2+\frac{1}{2\ell-1}}$, then a.a.s.
    $$
    \ex(G^r_{n,p}, C^r_{2\ell})\le cpn^{r-1}.
    $$

Let $X_m$ be the number of $C^r_{2\ell}$-free subgraphs of $G^r_{n,p}$ with $m$ edges and let $m=cpn^{r-1}$ for some constant $c$ sufficiently large. We have
$$
m\ge c\l(\log n\r)^{2r-2}n^{\frac{2\ell}{2\ell-1}}.
$$
Given that $K$ is sufficiently large, we can apply Theorem~\ref{theorem:container2_ge4} with $t=K$. This gives a constant $c_1>0$ and a collection $\C$ of subgraphs of $K^r_n$ such that
$$
|\C|\le \exp\l(c_1 (\log n)^{2r-2}n^{\frac{2\ell}{2\ell-1}}\r)\le\exp\l(\frac{c_1}{c}m\r);
$$
$e(G)\le Kn^{r-1}$ for each $G\in\C$;
every $C^r_{2\ell}$-free subgraph of $K^r_n$ is a subgraph of some $G\in\C$. In particular, every $m$-edge $C^r_{2\ell}$-free subgraph of $K^r_n$ is an $m$-edge subgraph of some $G\in\C$. So the expectation of $X_m$,
$$
\begin{aligned}
    \E[X_m]&\le |\C|\binom{Kn^{r-1}}{m}p^m\\
    &\le \exp \l(\frac{c_1}{c}m+\log\l(\frac{eKn^{r-1}}{m}\r)\cdot m+\log p\cdot m\r)\\
    &= \exp\l(\frac{c_1}{c}m+\log\l(\frac{eK}{cp}\r)\cdot m+\log p\cdot m\r)\\
    &= \exp\l(m\l(\frac{c_1}{c}+1+\log K-\log{c}\r)\r).
\end{aligned}
$$
Given that $c$ is sufficiently large, we have $\frac{c_1}{c}+1+\log K-\log{c}<0$. Then $\E[X_m]\to 0$ as $n\to\infty$. By Markov's Inequality,
$$
\P[X_m\ge 1]\le \E[X_m]\to 0,~\text{as~}n\to\infty.
$$
This implies that a.a.s.
$$
\ex(G^r_{n,p}, C^r_{2\ell})\le m=cpn^{r-1}.
$$
\end{proof}

% \begin{thmB}
%     For every $\ell\ge2$, there exists $c>0$ such that the following holds. Let
%     $$
%     p_0=cn^{-\frac{4\ell-3}{4\ell-2}}(\log n)^{\frac{4\ell^2+\ell-4}{2\ell-2}},~p_1=n^{-\frac{(\ell-1)(4\ell-3)}{4\ell^2-5\ell+2}}(\log n)^{3+\frac{2\ell-2}{4\ell^2-5\ell+2}}.
%     $$
%     Then a.a.s.
%     $$
%     \ex(G^3_{n,p}, C^3_{2\ell})\le 
%     \l\{
%     \begin{aligned}
%         &cp^{\frac{2(\ell-1))}{\ell(4\ell-3)}}n^{\frac{\ell+1}{\ell}}(\log n)^{\frac{2(2\ell-1)^2}{\ell(4\ell-3)}},~~&\text{if}~p_0\le p<p_1;\\
%         &cpn^2,~~&\text{if}~p\ge p_1.
%     \end{aligned}
%     \r.
%     $$
% \end{thmB}   

\begin{proof}[Proof of Theorem~\ref{theorem:main_3}]
Recall that we wish to prove that for every $\ell\ge2$, there exists $c>0$ such that the following holds. Let
    $$
    p_0=cn^{-\frac{4\ell-3}{4\ell-2}}(\log n)^{\frac{4\ell^2+\ell-4}{2\ell-2}},~p_1=n^{-\frac{(\ell-1)(4\ell-3)}{4\ell^2-5\ell+2}}(\log n)^{3+\frac{2\ell-2}{4\ell^2-5\ell+2}}.
    $$
    Then a.a.s.
    $$
    \ex(G^3_{n,p}, C^3_{2\ell})\le 
    \l\{
    \begin{aligned}
        &cp^{\frac{2(\ell-1)}{\ell(4\ell-3)}}n^{1+\frac{1}{\ell}}(\log n)^{3-\frac{4\ell-4}{\ell(4\ell-3)}},~~&\text{if}~p_0\le p<p_1;\\
        &cpn^2,~~&\text{if}~p\ge p_1.
    \end{aligned}
    \r.
    $$

Let $X_m$ be the number of $C^3_{2\ell}$-free subgraphs of $G^3_{n,p}$ with $m$ edges.

When $p\ge p_1$, let $m=cpn^2$ for some constant $c$ sufficiently large. We have
$$
m\ge cp_1n^2=c\l(\log n\r)^{3+\frac{2\ell-2}{4\ell^2-5\ell+2}}n^{\frac{4\ell^2-3\ell+1}{4\ell^2-5\ell+2}}.
$$
Given that $K$ is sufficiently large, we can apply Theorem~\ref{theorem:container2_3} with $t=K$. This gives a constant $c_1>0$ and a collection $\C_1$ of subgraphs of $K^3_n$ such that
$$
|\C_1|\le \exp\l(c_1 (\log n)^{3+\frac{2\ell-2}{4\ell^2-5\ell+2}}n^{\frac{4\ell^2-3\ell+1}{4\ell^2-5\ell+2}}K^{-\frac{2\ell-2}{4\ell^2-5\ell+2}}\r)\le\exp\l(\frac{c_1}{c}K^{-\frac{2\ell-2}{4\ell^2-5\ell+2}}m\r);
$$
$e(G)\le Kn^2$ for each $G\in\C_1$;
every $C^3_{2\ell}$-free subgraph of $K^3_n$ is a subgraph of some $G\in\C_1$. In particular, every $m$-edge $C^3_{2\ell}$-free subgraph of $K^3_n$ is an $m$-edge subgraph of some $G\in\C_1$. So the expectation of $X_m$,
$$
\begin{aligned}
    \E[X_m]&\le |\C_1|\binom{Kn^2}{m}p^m\\
    &\le \exp \l(\frac{c_1}{c}K^{-\frac{2\ell-2}{4\ell^2-5\ell+2}}m+\log\l(\frac{eKn^2}{m}\r)\cdot m+\log p\cdot m\r)\\
    &= \exp\l(\frac{c_1}{c}K^{-\frac{2\ell-2}{4\ell^2-5\ell+2}}m+\log\l(\frac{eK}{cp}\r)\cdot m+\log p\cdot m\r)\\
    &= \exp\l(m\l(\frac{c_1}{c}K^{-\frac{2\ell-2}{4\ell^2-5\ell+2}}+1+\log K-\log{c}\r)\r).
\end{aligned}
$$
Given that $c$ is sufficiently large, we have $\frac{c_1}{c}K^{-\frac{2\ell-2}{4\ell^2-5\ell+2}}+1+\log K-\log{c}<0$. Then $\E[X_m]\to 0$ as $n\to\infty$. By Markov's Inequality,
$$
\P[X_m\ge 1]\le \E[X_m]\to 0,~\text{as~}n\to\infty.
$$
This implies that a.a.s.
$$
\ex(G^3_{n,p}, C^3_{2\ell})\le m=cpn^2.
$$

When $p_0\le p< p_1$, let $m=cp^{\frac{2\ell-2}{\ell(4\ell-3)}}n^{\frac{\ell+1}{\ell}}(\log n)^{3-\frac{4\ell-4}{\ell(4\ell-3)}}$ for some sufficiently large constant $c$, and let
$$
t=p^{-\frac{4\ell^2-5\ell+2}{\ell(4\ell-3)}}n^{-\frac{\ell-1}{\ell}}(\log n)^{3-\frac{4\ell-4}{\ell(4\ell-3)}}K.
$$
Since $p_0\le p< p_1$, we have
$K<t\le n^{\frac{1}{2(2\ell-1)}}(\log n)^{-\frac{4\ell^2-7\ell+4}{2(\ell-1)}}c^{-\frac{4\ell^2-5\ell+2}{\ell(4\ell-3)}}K\le n^{\frac{1}{2(2\ell-1)}}(\log n)^{-\frac{4\ell^2-7\ell+4}{2(\ell-1)}}$ given that $c$ is sufficiently large.

By Theorem~\ref{theorem:container2_3} with $t$ above, we have a collection $\C_2$ of subgraphs of $K^3_n$ such that
$$
|\C_2|\le \exp\l(c_1 (\log n)^{3+\frac{2\ell-2}{4\ell^2-5\ell+2}}n^{\frac{4\ell^2-3\ell+1}{4\ell^2-5\ell+2}}t^{-\frac{2\ell-2}{4\ell^2-5\ell+2}}\r)=\exp\l(\frac{c_1}{c}K^{-\frac{2\ell-2}{4\ell^2-5\ell+2}}m\r);
$$
$e(G)\le tn^2$ for each $G\in\C_1$; every $C^3_{2\ell}$-free subgraph of $K^3_n$ is a subgraph of some $G\in\C_1$. In particular, every $m$-edge $C^3_{2\ell}$-free subgraph of $K^3_n$ is an $m$-edge subgraph of some $G\in\C_1$. So the expectation of $X_m$,
$$
\begin{aligned}
    \E[X_m]&\le |\C_2|\binom{tn^2}{m}p^m\\
    &\le \exp \l(\frac{c_1}{c}K^{-\frac{2\ell-2}{4\ell^2-5\ell+2}}m+\log\l(\frac{etn^2}{m}\r)\cdot m+\log p\cdot m\r)\\
    &= \exp\l(\frac{c_1}{c}K^{-\frac{2\ell-2}{4\ell^2-5\ell+2}}m+\log\l(\frac{eK}{cp}\r)\cdot m+\log p\cdot m\r)\\
    &= \exp\l(m\l(\frac{c_1}{c}K^{-\frac{2\ell-2}{4\ell^2-5\ell+2}}+1+\log K-\log{c}\r)\r)\to 0~\text{as~}n\to\infty.
\end{aligned}
$$
By Markov's Inequality,
$$
\P[X_m\ge 1]\le \E[X_m]\to 0,~\text{as~}n\to\infty.
$$
This implies that a.a.s.
$$
\ex(G^3_{n,p}, C^3_{2\ell})\le m= cp^{\frac{2\ell-2}{\ell(4\ell-3)}}n^{\frac{\ell+1}{\ell}}(\log n)^{3-\frac{4\ell-4}{\ell(4\ell-3)}}.
$$
\end{proof}

\section{Random Tur\'an Theorems for Berge even cycles}\label{section:berge}

In this section, we adapt the argument of Section~\ref{section:supersaturation}, 
\ref{section:container} and \ref{section:upper} to obtain a random Tur\'an theorem for Berge even cycles, that is, Theorem~\ref{Theorem:Berge}. Since the ideas are almost the same as those used for linear even cycles, we will only sketch the proof.

Observe that if we replace $C^3_{2\ell}$ in Theorem~\ref{theorem:supersaturation_3} with $\B^3_{2\ell}$ and replace the condition $t\ge K$ with $t\ge n^{-1+1/\ell+o(1)}$, then the theorem will still be true. This is because, when extending a $C_{2\ell}$ to a $C^3_{2\ell}$ we require the new vertices to be distinct; but when extending a $C_{2\ell}$ to a 3-uniform Berge $2\ell$-cycle we don't require the new vertices to be distinct. This means that we don't require all $2$-shadows in the selected ``large'' 2-shadow graph to have sufficiently large codegree. Hence, for Berge even cycles, we can skip the first paragraph after Case 2 in the proof of Theorem~\ref{theorem:supersaturation_3}. More precisely, by adapting the proof of Theorem~\ref{theorem:supersaturation_3} and Theorem~\ref{theorem:supersaturation_ge4}, we can obtain the following balanced supersaturation theorem for $\B^r_{2\ell}$. 

\begin{thm}\label{theorem:supersaturation_Berge}
    For every $r\ge 3$ and $\ell\ge 2$, there exist $K, c>0$ such that the following holds for all $n$ sufficiently large and $t\ge K(\log n)^{r-2}$. Given an $r$-graph $H$ with $n$ vertices and $tn^{1+1/\ell}$ edges, there exists a collection $\S$ of copies of elements in $\B^r_{2\ell}$ in $H$ such that, $,~\forall1\le j\le 2\ell$,
    $$
    \Delta_j(\S)\le\frac{c|S|(\log n)^{r-2}}{tn^{1+1/\ell}}\l(\max\l\{\l(\frac{t}{(\log n)^{r-2}}\r)^{-\frac{2(r-1)\ell^2-r\ell}{2(r-1)\ell^2-(r+2)\ell+2}},~\l(\frac{t}{(\log n)^{r-2}}\r)^{-1}n^{-\frac{\ell-1}{\ell(2\ell-1)}}\r\}\r)^{j-1}.
    $$
\end{thm}

\begin{proof}
    We use induction on $r$. When $r=2$, the statement is true by Theorem~\ref{theorem:cycle_supersaturation}. Now for $r\ge3$, suppose that the statement is true for $r-1$. We can assume that $H$ is $r$-partite with partition $V(H)=V_1\cup\dots V_r$. Let
    $$
    A=\l(\frac{t}{(\log n)^{r-2}}\r)^{\frac{\ell(2\ell-1)}{2(r-1)\ell^2-(r+2)\ell+2}}.
    $$
    We now partition $H$ into subgraphs $F$ and $F_{\tau,a}$, where $0\le a\le \log A$ and $\tau\in[r]^{r-1}$, as in the proof of Theorem~\ref{theorem:supersaturation_ge4}, such that every $(r-1)$-shadow of $F$ has codegree at least $A$ and every $(r-1)$-shadow in $\cup_{i\in\tau}V_i$ of $F_{\tau,a}$ has codegree at least $2^a$ and at most $2^{a+1}-1$. The proof splits into 2 cases.

    \noindent\textbf{Case 1:} $e(F)\ge e(H)/(\log n)^{r-2}$.\\

    In this case, we can use ``greedy expansion'' to find a collection of Berge $2\ell$-cycles. Precisely, applying Lemma~\ref{lem: expansion} on $F$ gives us a collection $\S$ of copies of $C^3_{2\ell}$ in $H$ such that, for $1\le j\le 2\ell-1$,
    $$
    \Delta_j(\S)\le O\l(\frac{|\S|}{e(F)}\l(A^{-r+1}\r)^{j-1}\r).
    $$
    Further, since $\Delta_{2\ell}(\S)=1$, we have, for $1\le j\le 2\ell$,
    $$
    \Delta_j(\S)\le O\l(\frac{|\S|}{e(F)}\l(A^{-\frac{(2\ell-1)r-2\ell}{2\ell-1}}\r)^{j-1}\r)\le O\l(\frac{|\S|(\log n)^{r-2}}{tn^{1+1/\ell}}\l(\frac{t}{(\log n)^{r-2}}\r)^{-\frac{2(r-1)\ell^2-r\ell}{2(r-1)\ell^2-(r+2)\ell+2}\cdot(j-1)}\r).
    $$

\medskip

    \noindent\textbf{Case 2:} $e(F)< e(H)/(\log n)^{r-2}$.\\

    In this case, by the Pigeonhole Principle, there exists a pair $(\tau,a)$ such that  $e(F_{\tau,a})\ge\Omega\l(\frac{tn^{1+1/\ell}}{\log n}\r)$. Let $F'=F_{\tau, a}$ and let $D=2^a$. Further, let $G$ be the $(r-1)$-graph on $\cup_{i\in\tau}V_i$ consisting of all $(r-1)$-shadows of $F'$ in $\cup_{i\in\tau}V_i$. Then $e(G)\ge e(F')/(2D)\ge \Omega\l(tn^{1+1/\ell}/(D\log n)\r)$. Let $t'=e(G)/n^{1+1/\ell}$. Then we have 
    \begin{equation}\label{equation:Berge_t'}
        t'\ge \Omega\l(\frac{t}{D\log n}\r).
    \end{equation}
    Further, by using $D\le A$ and $t\ge \Omega((\log n)^{r-2})$, one can check that $t'\ge\Omega \l((\log n)^{r-3}\r)$. Hence by the inductive hypothesis, we can find a balanced collection $\S'$ of copies of $(r-1)$-uniform Berge $2\ell$-cycles in $G$ such that, $\forall 1\le j\le 2\ell$,
    \begin{equation}\label{equation:Berge_shadow}
        \Delta_j(\S')\le O\l(\frac{|S'|(\log n)^{r-3}}{t'n^{1+1/\ell}}\l(\max\l\{\l(\frac{t'}{(\log n)^{r-3}}\r)^{-\frac{2(r-2)\ell^2-(r-1)\ell}{2(r-2)\ell^2-(r+1)\ell+2}},~\l(\frac{t'}{(\log n)^{r-3}}\r)^{-1}n^{-\frac{\ell-1}{\ell(2\ell-1)}}\r\}\r)^{j-1}\r).
    \end{equation}

    Since each $(r-1)$ edge of $G$ is contained in at least $D$ $r$-edges of $F'$, we can extend each $(r-1)$-uniform Berge $2\ell$-cycle in $G$ into at least $D^{2\ell}$ copies of $r$-uniform Berge $2\ell$-cycles in $H$, which gives us a collection $\S$ of Berge $2\ell$-cycles such that
    \begin{equation}\label{equation:Berge_shadow_expansion}
        |\S|\ge\Omega\l(|\S'|D^{2\ell}\r).
    \end{equation}
    Further, since each $(r-1)$ edge of $G$ is contained in at most $2D$ $r$-edges of $F'$, together with (\ref{equation:Berge_t'}), (\ref{equation:Berge_shadow}) and (\ref{equation:Berge_shadow_expansion}), we know that, $\forall 1\le j\le 2\ell$,
    $$
    \begin{aligned}
    &\Delta_j(\S)\le O\l(\Delta_j(\S')D^{2\ell-j}\r)\\
    &\le O\l(\frac{|\S|(\log n)^{r-2}}{tn^{1+1/\ell}}\l(\max\l\{\l(\frac{D^{-\frac{2\ell-2}{2(r-2)\ell^2-(r-1)\ell}}t}{(\log n)^{r-2}}\r)^{-\frac{2(r-2)\ell^2-(r-1)\ell}{2(r-2)\ell^2-(r+1)\ell+2}},~\l(\frac{t}{(\log n)^{r-2}}\r)^{-1}n^{-\frac{\ell-1}{\ell(2\ell-1)}}\r\}\r)^{j-1}\r).
    \end{aligned}
    $$
    Finally, note that
    $$
    D\le A=\l(\frac{t}{(\log n)^{r-2}}\r)^{\frac{\ell(2\ell-1)}{2(r-1)\ell^2-(r+2)\ell+2}},
    $$
    hence we have, $\forall 1\le j\le 2\ell$,
    $$
    \Delta_j(\S)\le\l(\frac{|S|(\log n)^{r-2}}{tn^{1+1/\ell}}\l(\max\l\{\l(\frac{t}{(\log n)^{r-2}}\r)^{-\frac{2(r-1)\ell^2-r\ell}{2(r-1)\ell^2-(r+2)\ell+2}},~\l(\frac{t}{(\log n)^{r-2}}\r)^{-1}n^{-\frac{\ell-1}{\ell(2\ell-1)}}\r\}\r)^{j-1}\r).
    $$
\end{proof}

Using Theorem~\ref{theorem:supersaturation_Berge} and Theorem~\ref{theorem:container} with coefficients 
$$
    L=\frac{ctn^{1+1/\ell}}{(\log n)^{r-2}},~\text{and}~B=(\log n)^{r-2}\max\l\{\l(\frac{t}{(\log n)^{r-2}}\r)^{-\frac{2\ell-2}{2(r-1)\ell^2-(r+2)\ell+2}}n^{1+1/\ell},n^{1+\frac{1}{2\ell-1}}\r\},
$$
we can obtain the following container theorems for $\B^r_{2\ell}$. We omit the details.

\begin{thm}\label{theorem:container1_Berge}
    For every $\ell\ge2$ and $r\ge3$, there exist $\epsilon, K>0$ such that the following holds for all sufficiently large $n$ and every $t\ge K(\log n)^{r-2}$. Given an $r$-graph $H$ with $n$ vertices and $tn^{1+1/\ell}$ edges there exists a collection $\C$ of at most 
    $$
    \exp\l(\frac{(\log n)^{r-1}}{\epsilon}\max\l\{\l(\frac{t}{(\log n)^{r-2}}\r)^{-\frac{2\ell-2}{2(r-1)\ell^2-(r+2)\ell+2}}n^{1+1/\ell},n^{1+\frac{1}{2\ell-1}}\r\}\r)
    $$
    subgraphs of $H$ such that
    \begin{itemize}
        \item[(a)] every $\B^r_{2\ell}$-free subgraph of $H$ is a subgraph of some $G\in\C$;
        \item[(b)] $e(G)\le \l(1-\frac{\epsilon}{(\log n)^{r-2}}\r)tn^{1+1/\ell}$, for each $G\in\C$. 
    \end{itemize}
\end{thm}

By using Theorem~\ref{theorem:container1_Berge} repeatedly, we obtain the following theorem.

\begin{thm}\label{theorem:container2_Berge}
    For every $\ell\ge2$ and $r\ge3$, there exist $c, K>0$ such that for all sufficiently large $n$ and every $t\ge K(\log n)^{r-2}$, there exists a collection $\C$ of $r$-graphs on $[n]$ such that 
    \begin{enumerate}
        \item[(a)]
        $$
        |\C|\le \exp\l(c{(\log n)^{2r-3}}\max\l\{\l(\frac{t}{(\log n)^{r-2}}\r)^{-\frac{2\ell-2}{2(r-1)\ell^2-(r+2)\ell+2}}n^{1+1/\ell},\log n\cdot n^{1+\frac{1}{2\ell-1}}\r\}\r);
        $$
        \item[(b)] $e(G)\le tn^{1+1/\ell}$, $\forall G\in\C$;
        \item[(c)] every $\B^r_{2\ell}$-free $r$-graph on $[n]$ is a subgraph of some $G\in\C$.
    \end{enumerate}
\end{thm}

Now Theorem~\ref{Theorem:Berge} follows easily from Theorem~\ref{theorem:container2_Berge} using the method in Section~\ref{section:upper}. We omit the details.

\section{Concluding remarks}
\begin{itemize}
    \item Given that the proof of Theorem~\ref{theorem:main_3} is based on Theorem~\ref{theorem:cycle_supersaturation}, the supersaturation of $C_{2\ell}$ by Morris and Saxton~\cite{MORRIS2016534}, it is natural to ask if it is possible to improve Theorem~\ref{theorem:cycle_supersaturation} in order to close the gap in Theorem~\ref{theorem:main_3}. Unfortunately, Theorem~\ref{theorem:cycle_supersaturation} is essentially tight conditional on Conjecture~\ref{conjecture:Erdos-Simonovits} which is widely believed to be true. Hence, in order to improve Theorem~\ref{theorem:main_3}, one might need some new ideas.
    \item Given an $r$-graph $H$, let $N(n,m,H)$ be the number of $H$-free $r$-graph on $[n]$ with $m$ edges. In the proof of Theorem~\ref{theorem:main_ge4}, Theorem~\ref{theorem:main_3} and Theorem~\ref{Theorem:Berge}, we in fact have also implicitly proved the following theorems, which may be of independent interest.
    \begin{thm}
        For $r\ge 4$ and $\ell\ge2$, when $n^{1+\frac{1}{2\ell-1}+o(1)}\le m\le o(n^{r-1})$,
        \begin{equation}\label{equation:Nm_ge4}
            N(n,m,C^r_{2\ell})\le \l(\frac{n^{r-1}}{m}\r)^{(1+o(1))m}.
        \end{equation}
    \end{thm}
    Note that (\ref{equation:Nm_ge4}) is tight up to the $o(1)$ in the exponent, since otherwise a better upper bound would improve Theorem~\ref{theorem:main_ge4}, which we know is tight.
    \begin{thm}
        For $\ell\ge 2$, when $n^{\frac{4\ell^2-3\ell+1}{4\ell^2-5\ell+2}+o(1)}\le m\le o(n^2)$,
        \begin{equation}\label{equation:Nm_3_1}
            N(n,m,C^3_{2\ell})\le \l(\frac{n^{2}}{m}\r)^{(1+o(1))m};
        \end{equation}
        when $n^{1+\frac{1}{2\ell-1}+o(1)}\le m\le n^{\frac{4\ell^2-3\ell+1}{4\ell^2-5\ell+2}-o(1)}$,
        \begin{equation}\label{equation:Nm_3_2}
            N(n,m,C^3_{2\ell})\le \l(\frac{n^{\ell+1}}{m^{\ell}}\r)^{\l(\frac{4\ell-3}{2\ell-2}+o(1)\r)m}.
        \end{equation}
    \end{thm}
    Similarly, (\ref{equation:Nm_3_1}) is tight up to the $o(1)$ term; but (\ref{equation:Nm_3_2}) is not. In order to close the gap for the random Tur\'an number of $C^3_{2\ell}$, we expect to show that (\ref{equation:Nm_3_1}) holds for all $n^{1+\frac{1}{2\ell-1}+o(1)}\le m\le o(n^2)$.
    \begin{thm}
    For $r\ge 3$ and $\ell\ge2$, when $n^{1+\frac{1}{2\ell-1}+o(1)}\le m\le n^{1+\frac{1}{\ell}-o(1)}$, \
    $$
    N\l(n,m,\B^r_{2\ell}\r)\le \l(\frac{n^{\ell+1}}{m^\ell}\r)^{\l(\frac{2(r-1)\ell-r}{2\ell-2}+o(1)\r)m}.
    $$
    \end{thm}
    This theorem is not tight, since the random Turán upper bounds (Theorem~\ref{Theorem:Berge}) that it implies are not tight, as commented after Theorem~\ref{Theorem:Berge}.
    \item In this paper, we adopt the strategy of induction and codegree dichotomy to prove new balanced supersaturation theorems for $C^{r}_{2\ell}$. This method can also be applied to some other families of hypergraphs; for example, the family of odd linear cycles, the family of expansions of tight trees, and the family of expansion of complete hypergraphs. In fact, it is possible to prove a general theorem which gives reasonably good upper bounds for all hypergraphs mentioned above. We do those in a followup paper~\cite{nie2023random}.
\end{itemize}

\section{Acknowledgements}
The author would like to express gratitude to Jacques Verstra\"ete for introducing him to the random Tur\'an problems. Additionally, the author wishes to thank Dhruv Mubayi, Hehui Wu, and Liana Yepremyan for helpful discussions. Special thanks are owed to Sam Spiro for bringing to attention that the proof for linear even cycles could be readily modified to accommodate Berge even cycles, resulting in Theorem~\ref{Theorem:Berge}. Finally, the author would like to thank the anonymous referees for their thorough reading of the preprint and their helpful comments and suggestions.

\bibliographystyle{abbrv}
\bibliography{refs}

\end{document}